\documentclass[10pt]{amsart}

\allowdisplaybreaks[4]
\usepackage{amsmath,amsthm,amscd,amssymb,amsfonts, amsbsy}
\usepackage{latexsym, color}
\usepackage{pxfonts}
\usepackage{mathrsfs}
\IfFileExists{bbm.sty}{\usepackage{bbm}}{}
\usepackage{mathtools}

\usepackage{lineno,hyperref}
\modulolinenumbers[5]

\usepackage{enumitem}

\usepackage[final,notref,notcite]{showkeys}
\usepackage{todonotes}

\theoremstyle{plain}
\newtheorem{theorem}{Theorem}
\newtheorem{lemma}{Lemma}
\newtheorem{corollary}{Corollary}
\newtheorem{proposition}{Proposition}

\theoremstyle{definition}

\theoremstyle{remark}
\newtheorem{remark}{Remark}

\newcommand{\dv}{\operatorname{div}}
\newcommand{\divg}{\operatorname{div}}

\newcommand{\supp}{\operatorname{supp}}

\numberwithin{equation}{section}

\newcommand{\bN}{\mathbb{N}}

\newcommand{\bR}{\mathbb{R}}

\newcommand{\vphi}{\varphi}

\newcommand{\ep}{\epsilon}
\newcommand{\lt}{\left}
\newcommand{\rt}{\right}

\newcommand\cC{\mathcal{C}}

\newcommand\cE{\mathcal{E}}

\newcommand\cI{\mathcal{I}}
\newcommand\cJ{\mathcal{J}}

\newcommand\cM{\mathcal{M}}

\newcommand\sF{\mathscr{F}}

\newcommand\te{\theta}

\newcommand{\ud}{\textnormal{d}}
\newcommand{\dt}{\ud t}

\newcommand{\dtau}{\ud \tau}

\definecolor{MutedGreen}{HTML}{397052}

\newcommand{\Del}{\Delta}

\newcommand{\Lam}{\Lambda}
\newcommand{\Ome}{\Ome}

\providecommand{\norm}[1]{\lVert#1\rVert}
\providecommand{\Norm}[1]{\left\lVert#1\right\rVert}

\newcommand{\p}{\partial}

\newcommand{\qand}{\quad \text{and} \quad}

\begin{document}
\title[Stability and temporal decay for the tropical climate model]
{Higher-order Sobolev stability and temporal decay estimate\\ for the tropical climate model with small $H^1$ data}
\date{November 7, 2026}

\author[H. In]{Hyunjin In}
\address[Hyunjin In]{Research Institute of Basic Sciences, Ajou University, Suwon 16499, Republic of Korea}
\email{hyunjinin@ajou.ac.kr}

\author[D. Kim]{Dong-ha Kim}
\address[Dong-ha Kim]
{Research Institute of Mathematics, Seoul National University, 1 Gwanak-ro, Gwanak-gu, Seoul 08826, 
Republic of Korea,}
\email{kdh.pde@snu.ac.kr}

\author[J. Kim]{Junha Kim}
\address[Junha Kim]{Department of mathematics, Ajou University 206, World cup-ro, Yeongtong-gu, Suwon-si, Gyeonggi-do, Republic of Korea}
\email{junha02@ajou.ac.kr}


\begin{abstract}
We study the stability and temporal decay of solutions to the two-dimensional tropical climate model with Laplacian dissipation in both velocity equations and no thermal diffusion.
For initial data in $H^s (\mathbb{R}^2)$ ($s \ge 1$) with sufficiently small $H^1$ norm, we establish a unique global solution satisfying uniform $H^s$ bounds and time-integrated estimates for the velocity and temperature.
The smallness is independent of $s$, and the higher Sobolev norms of the initial data need not be small.
Furthermore, without additional spatial integrability or negative-order regularity of the initial data, we prove that the $\dot{H}^s (\mathbb{R}^2)$ norms of the solutions are $o(t^{- s / 2})$ as $t \to \infty$.
We also show that this decay rate is sharp: no additional decay factor tending to zero works for all such initial data.
\end{abstract}

\maketitle


\section{Introduction and main results}
\label{S1}

\subsection{Introduction}
Frierson--Majda--Pauluis \cite{MR2119930} introduced the tropical climate model on $\bR^2$ to study the interaction between large-scale
tropical circulation and precipitation:
\begin{align*}
\begin{cases}
    u_t + (u \cdot \nabla) u + \nabla p 
    = - \dv (v \otimes v), \\
    v_t + (u \cdot \nabla) v + (v \cdot \nabla) u = \nabla \theta, \\
    \theta_t + u \cdot \nabla \theta = \dv v, \\
    \dv u =0.
\end{cases}
\end{align*} 
Here, the vector fields $u$ and $v$ denote the barotropic and the first baroclinic modes of the velocity field (or simply, baroclinic modes), respectively. 
The scalar functions $\theta$ and $p$ denote the temperature and pressure.

The dynamical core of the model is obtained by a vertical Galerkin truncation of the hydrostatic Boussinesq equations, retaining the barotropic and baroclinic modes.
In this truncation, the baroclinic velocity profile changes sign between the lower and upper troposphere, representing oppositely directed horizontal motions;
see \cite[Section~1.4]{Khouider2019}.
Retaining both velocity modes preserves the nonlinear interaction
between the vertically averaged flow and its baroclinic
component.
The baroclinic motion contributes a momentum flux to the barotropic
equation, while the barotropic velocity enters the baroclinic
equation through advection and velocity-gradient coupling.
This interaction provides a mechanism for momentum transfer
between the modes and is relevant to models of
tropical—extratropical interactions
\cite[Section~2.2]{MR2119930}.
A derivation of the viscous velocity—temperature equations is given in \cite[Section~6]{LiTiti2018}.

In this work, we consider the tropical climate model with Laplacian diffusion in both the barotropic mode $u$ and the baroclinic mode $v$, but no thermal diffusion in $\bR^2$:
\begin{align}
\label{E11}
\begin{cases}
    u_t + (u \cdot \nabla) u + \nabla p  
   -  \Delta u
    = - \dv (v \otimes v), \\
    v_t + (u \cdot \nabla) v + (v \cdot \nabla) u 
   -  \Delta v
    = \nabla \theta, \\
    \theta_t + u \cdot \nabla \theta 
    = \dv v, \\
    \dv u =0.
\end{cases}
\end{align} 

The tropical climate model has attracted considerable attention in the mathematical literature.
In a seminal work \cite{MR3479523}, Li--Titi first established global existence and uniqueness of strong solutions to the system \eqref{E11}, for initial data in $H^1 (\bR^2)$ satisfying $\dv u_0 = 0$ and
$\theta_0 \in L^\infty (\bR^2)$.
A key ingredient in their proof is the pseudo-baroclinic velocity $w:=v-\nabla(-\Delta)^{-1}\theta$, which is a combination of $v$ and $\theta$, and allows them to overcome the difficulty from the absence
of thermal diffusion.
Liu--Zhao \cite{MR4969841} subsequently removed the assumption $\te_0 \in L^{\infty}(\bR^2)$ by proving uniqueness in Lagrangian coordinates.
For higher Sobolev regularity, Ma--Jiang--Wan \cite{MaJiangWan2016} established local existence and uniqueness for \eqref{E11} in $H^s (\bR^2)$ with $s > 1$. 
Dong--Wu--Ye \cite{DongWuYe2019} established global existence and uniqueness in $H^s (\bR^2)$ with $s > 2$.

Next, we consider the stability problem of whether suitable smallness assumptions on the initial data ensure uniform-in-time bounds for the solutions.
In the Sobolev space setting, the constant-coefficient case considered by Ye--Zhu
\cite{YeZhu2020} covers \eqref{E11}.
For initial data in $H^s(\bR^2)$ with $s>1$, they established
global-in-time regularity result provided that
$\|v_0\|_{H^s}+\|\theta_0\|_{H^s}$ is sufficiently small depending on $\|u_0\|_{H^s}$.
In the Besov space setting, Li--Zhai--Yin \cite{MR3973083} obtained a unique global solution in low regularity spaces when $\|u_0\|_{\dot B_{2,1}^0} +\|v_0\|_{\dot B_{2,1}^0}+\|\te_0\|_{\dot B_{2,1}^0\cap \dot B_{2,1}^1}$ is sufficiently small; see also \cite{ZhaiChen2020}.

For tropical climate models with fractional or partial dissipation,
we refer the reader to
\cite{Ye2017,DongWangWuZhang2019,DongWangWuYeZhang2019,DongWuYe2020,Ye2020,LiYeZhong2026}
and the references therein.
For models with linear or nonlinear damping, see \cite{Wan2016,MaWan2018,LiYu2020,ChenYuan2026}
and the references therein.
For models with temperature-dependent coefficients, see \cite{YeZhu2019,DongLiXuYe2021,LiXuYe2022,MR5029968} and the references therein.
For decay estimates in diffusive models, stability near Couette flow, and other works, see \cite{LiXiao2019,LiDengShang2021,NiuWu2022,NiuWu2025,RenYin2026,NiuWuWuXu2026,PaicuZhu2021,YangHanZhaoCaraballo2023,DingWu2023} and the references therein.

\subsection{Main theorem}
In this paper, we establish stability and temporal decay estimate for the solution to \eqref{E11} in $H^s(\bR^2)$ when $s\geq1$, 
under a smallness assumption only on the initial $H^1$ norm. Moreover, we explain the sharpness of the decay estimate.

Throughout this paper, for any normed space $X$, we use the shorthand
$\|(f,g,h)\|_X^2:=\|f\|_X^2+\|g\|_X^2+\|h\|_X^2$,
and similarly for pairs.

\begin{theorem}\label{Thm1}
    There exists a small absolute constant $\ep_0 > 0$ which satisfies the following: let $s\geq1$ and let the initial data $(u_0,v_0,\theta_0)$ belongs to $H^s(\bR^2)$ with $\dv u_0=0$. Suppose that
\[
 \|(u_0,v_0,\te_0)\|_{H^1(\bR^2)}<\ep_0.
\]

\begin{enumerate}[label=(\roman*)]
    \item 
Then, there exists a unique solution $(u, v, \te) \in \mathcal{C} ([0, \infty); H^s (\bR^2))$ which satisfies the following bound:
\begin{equation}
\label{E10}
\begin{aligned}
\sup_{t \in [0, \infty)}\|(u, v, \te)(t)\|_{H^s} 
&+ \left( \int_0^{\infty} 
\|(\nabla u , \nabla v )(t)\|_{H^s}^2
+ \|\nabla \te (t)\|_{H^{s-1}}^2 \dt \right)^{\frac{1}{2} }
\le C \|(u_0,v_0,\te_0)\|_{H^s (\bR^2)}.
\end{aligned}
\end{equation}
    \item
Moreover, there exists a non-increasing function $\cM_s:[0, \infty)\to[0, \infty)$ 
such that
\begin{equation}
\label{E13}
(1 + t)^{\frac{s}{2}}
\left(
\|(\Lam^s u,
\Lam^s v,
\Lam^s \te)(t)\|_{L^2}^2
+\int_t^{\infty}
\|(\Lam^{s+1}u
,\Lam^{s+1}v)(\tau)\|_{L^2}^2
\,
\dtau
\right)
^\frac{1}{2}
\le \cM_s(t),
\qquad \text{for all } t \ge 0,
\end{equation}
with
\begin{equation*}
    \sup_{t\in [0, \infty)}\cM_s(t) \leq C\|(u_0,v_0,\te_0)\|_{H^{s}(\bR^2)},\qquad\textrm{ and }\qquad \lim_{t\to\infty}\cM_s(t)=0. 
\end{equation*}
\end{enumerate}
Here, the constants $C$ depend only on $s$.
\end{theorem}
\begin{remark}
For $s>1$, the higher Sobolev norms of the initial data $\|(u_0,v_0,\te_0)\|_{H^s}$ need not be small.
\end{remark}

The decay rate in \eqref{E13} cannot be improved by any additional factor tending to zero, as the following theorem shows.
\begin{theorem}[Sharpness of the temporal decay]
\label{Thm2:sharp}
Let $\ep_0$ be as in Theorem~\ref{Thm1}.
For every $0 < \ep < \ep_0$ and every function
$h : [1, \infty) \to (0, \infty)$ with $\lim_{t \to \infty} h (t) = 0$,
there exists a radial function
$\te_0 \in \bigcap_{m = 1}^{\infty} H^m (\bR^2)$ satisfying
$0 < \|\te_0\|_{H^1} < \ep$ such that the solution of \eqref{E11}
with initial data $(0, 0, \te_0)$ satisfies $u \equiv 0$ and
\[
\limsup_{t \to \infty}
\frac{(1 + t)^{\frac{s}{2}}}{h (t)}
\left (
\|(\Lambda^s u, \Lambda^s v, \Lambda^s \te) (t)\|_{L^2}^2
+ \int_t^\infty
\|(\Lambda^{s + 1} u, \Lambda^{s + 1} v) (\tau)\|_{L^2}^2 d\tau
\right)^{\frac{1}{2}}
= \infty
\]
for every real $s \ge 1$.
The initial data may depend on $\ep$ and $h$, but are independent of $s$.
\end{theorem}

To briefly illustrate our proof strategy, consider the simplified system
\begin{equation}\label{toy_exam1}
\begin{cases}
    u_t -\Delta u=0, \qquad \dv u=0, \\
    v_t-\Delta v=\nabla\theta, \\
    \theta_t+u\cdot \nabla \te=\dv v,\\
\end{cases}
\end{equation}
We also assume that $\nabla u \in L^{1}_{\text{loc}}([0, \infty);L^\infty)$.
The $\dot{H}^1$ energy estimate reads
\begin{equation}\label{toy_estimate1}
    \frac{1}{2}\frac{\ud}{\ud t}
        \|(\nabla u, \nabla v, \nabla \te)\|_{L^2}^2
    +\|(\Delta u, \Delta v)\|_{L^2}^2\leq \|\nabla u\|_{L^\infty}\|\nabla\theta\|_{L^2}^2.    
\end{equation}
To recover the dissipation for $\theta$, we exploit its coupling
with $v$ by introducing a \textit{cross term}, as in
\cite{YeZhu2020,MR5029968}.
After some manipulation of \eqref{toy_exam1}, integration by parts and Young's inequality give
\[
    \frac{1}{2}\|\nabla\theta\|_{L^2}^2 - \frac{\ud}{\ud t}
    \int_{\bR^2}v\cdot\nabla\theta\,\ud x
    \leq
    \|\dv v\|_{L^2}^2 + \frac{1}{2}\|\Del v\|_{L^2}^2
+ \left|\int_{\bR^2}(u\cdot\nabla\theta)\dv v\,\ud x\right|.
\]
Combining this cross term with the $H^1$ energy therefore yields
a modified energy
\[
    \mathcal E(t):=
    \|(u(t), v(t), \te(t))\|_{H^1}^2
    -\frac{1}{2}\int_{\bR^2}v(t)\cdot\nabla\theta(t)\,\ud x,
\]
which is equivalent to the sum of the squared $H^1$ norms and satisfies
\[
    \frac{\ud}{\ud t}\mathcal E
    +\frac{1}{4}\left(
\|(\nabla u, \nabla v)\|_{H^1}^2
        +\|\nabla\theta\|_{L^2}^2
    \right)
\leq 2\|\nabla u \|_{L^{\infty}}\|\nabla \te\|_{L^2}^2 +\frac{1}{2}\left|\int_{\bR^2}(u\cdot\nabla\theta)\dv v\,\ud x\right|.
\]
This produces a non-negative $\norm{\nabla \theta}_{L^2}^2$ dissipation term. 
Under a smallness condition, it is easy to control the last term on the right-hand side.
In the space $H^s(\bR^2)$ with $s>1$, we can derive a similar $H^s$-estimate and enjoy the embedding $\|\nabla u\|_{L^{\infty}} \lesssim \|\nabla u\|_{H^{s}}$ to close the stability analysis, as in \cite{YeZhu2020,MR5029968}. 
However, the same approach does not work at $s=1$, since this embedding fails. Not surprisingly, the global well-posedness results in \cite{MR3479523,MR4969841} do not by themselves provide bounds for $\|\nabla u\|_{L^1([0,T];L^{\infty})}$ which is uniform in $T$; the available bound may grow with $T$ (see the estimate in \cite[Corollary~1]{MR3479523}).

To overcome this difficulty for the full system \eqref{E11}, we estimate the Duhamel formula for $u$ in Fourier space and obtain bounds for
$\|\nabla u\|_{L^p(0,T;L^\infty)}$ with $1<p<2$.
Combined with the dissipation $\|\nabla \te\|_{L^2}^2$ recovered from the cross term, these estimates allow us to close the $H^1$ stability estimate under a smallness assumption on the initial $H^1$ norm. 
See Lemma~\ref{Prop1} for details.

We next turn to temporal decay, namely, estimates of the form \eqref{E13}.  We restrict ourselves to the case $s=1$. To this end, we use a modified energy inequality, since \eqref{toy_estimate1} alone does not capture the mechanism of temporal decay.
Indeed, the same inequality holds for the fully decoupled system
\[
\begin{cases}
    u_t -\Delta u=0, \qquad \dv u=0, \\
    v_t-\Delta v=0, \\
    \theta_t=0,
\end{cases}
\]
in which $\te(t,\cdot)=\te_0$ regardless of how small its initial norm is. To establish decay, we therefore combine the dissipation recovered from the cross term with the stability estimates. 

The key step is a \textit{minimum-selection argument}.
For each sufficiently large $T$, choose
$r^*\in[T/2,T]$ minimizing the continuous map
$t \mapsto\|(\nabla u,\nabla v,\nabla\te)(t)\|_{L^2}$.
The time-integrated gradient bounds of $\norm{\nabla \theta}_{L^2}^2$ established in the stability analysis (see \eqref{eq1016wed}) give
\begin{equation*}
\begin{aligned}
\|(\nabla u, \nabla v,\nabla \te) (T)\|_{L^2}^2 
+\int_{r^*}^{T} \|(\Delta u, \Delta v) (\tau)\|_{L^2}^2\dtau &\leq C\|(\nabla u, \nabla v,\nabla \te) (r^*)\|_{L^2}^2 \\ 
&
\leq \frac{C}{T}\int_{T/2}^{T} \|(\nabla u, \nabla v,\nabla \te) (t)\|_{L^2}^2 \,
\dt\leq \frac{C}{T}\int_{T/2}^{\infty} \|(\nabla u, \nabla v,\nabla \te) (t)\|_{L^2}^2 \,\dt.
\end{aligned}
\end{equation*}
Since the last integral tends to zero, we obtain
$\|(\nabla u, \nabla v, \nabla \te)(t)\|_{L^2}=o(t^{-1/2})$. See the proof of Proposition~\ref{prop_stab_H1} for the details.
We now briefly discuss the sharpness of the decay estimates.
For two-dimensional Navier--Stokes flows with radially symmetric assumption, radial symmetry and incompressibility imply that the nonlinear term can be absorbed into the pressure; see \cite[p.~224]{Brandolese2004}. The symmetry-preserving structure of the system~\eqref{E11} yields a related reduction: if $u_0=v_0=0$ and $\theta_0$ is radially symmetric, then $u$ remains identically zero, and the remaining equations form a linear system; see Lemma~\ref{lem:S4} and the proof of Theorem~\ref{Thm2:sharp}. We can then prove that the decay estimates are sharp by a direct calculation.

These are the main ideas of the proofs.

\subsection{Outline of the paper}
This paper is organized as follows.
In Section~\ref{S2}, we prove Proposition~\ref{prop_stab_H1},
which establishes the stability and temporal decay estimates
of Theorem~\ref{Thm1} in the case $s=1$.
We also derive time-integral estimates on dyadic intervals
that will be used in the higher-order analysis.
In Section~\ref{S3}, we extend these results to $s>1$,
thereby completing the proof of Theorem~\ref{Thm1}.
Section~\ref{S4} proves Theorem~\ref{Thm2:sharp}, showing that the decay in \eqref{E13} cannot be improved by any additional factor tending to zero.

\subsection{Acknowledgement}
Hyunjin In was supported by the National Research Foundation of Korea(NRF) grant funded by the Korea government(MSIT) (RS-2024-00360798 and RS-2026-25566720) and Basic Science Research Program through the  National Research Foundation of Korea(NRF) funded by the Ministry of Education(RS-2021-NR060141). Dong-ha Kim was supported by the SSEN-Muju Fellow Research Grant and the Research Institute of Mathematics, Seoul National University. Junha Kim was supported by the National Research Foundation of Korea(NRF) grant funded by the Korea government(MSIT) (RS-2024-00360798).

\subsection{Disclosure on AI assistance}
The authors used GPT-6.1 during the research process, including for literature searches and exploratory
discussions of possible proof strategies.

The authors developed the strategy for the stability part of
Theorem~\ref{Thm1}, drawing on \cite{MR5029968,JKL2022}. 
In this part, all of the main ideas were generated by the authors, and all mathematical statements and proofs have been checked by the authors.
For the temporal decay estimates, the authors initially obtained the stated decay in $\dot H^\gamma$ only for $0<\gamma\leq s-1$ (see \cite{MR5029968}).
During subsequent discussions with GPT-6.1,
a \textit{minimum-selection argument} was suggested.
This suggestion was essential for removing the restriction on $\gamma$ and improving the bound for $\cM_s$. GPT-6.1 also suggested an idea for establishing the optimality of the temporal decay estimate.

The authors formulated the final theorem statements, wrote the proofs, independently verified all mathematical arguments and cited references, and take full responsibility for the contents of the paper.

\section{The special case: $s=1$}
\label{S2}
We introduce the following notation for use throughout the paper.
\begin{itemize}
\item
We denote $A \lesssim B$ if there exists a generic positive constant $C$ such that $A \le C B$.
To stress the dependence of $C$ on a parameter $M$, we write $A \lesssim_M B$.
\item
The letter $C$ represents a generic constant, which may change from
line to line.
\item
We simply write $\norm{f}_{L^p} = \norm{f}_{L^p ( \bR^2)}$ for $p \in [1,\infty]$,  and $\norm{f}_{H^s} = \norm{f}_{H^s ( \bR^2)}$ for $s >0$. 
\item
For a spatial normed space $X$ and $1\leq p<\infty$, we write
$\|f\|_{L^p(r,t;X)}:=(\int_r^t\|f(\tau)\|_X^p\,\ud\tau)^{1/p}$.
For $p=\infty$, we set
$\|f\|_{L^\infty(r,t;X)}:=\operatorname*{ess\,sup}_{\tau\in(r,t)}\|f(\tau)\|_X$.
\item
The Fourier transform of $f$ is denoted by $\sF f$.
\item  
For vector fields $a,b$, we write $(a \otimes b)$ as $(a_i b_j)_{1 \le i,j \le 2}$.
We denote $A : B = \sum_{i,j=1}^2(a_{ij} b_{ij})$ for $2 \times 2$ matrices $A = (a_{ij})$ and $B = (b_{ij})$.
\end{itemize}

To avoid repetition, we first recall the global-in-time existence result for $H^1$ strong solutions from \cite[Theorem 1.1]{MR3479523} and \cite[Theorem 1.1]{MR4969841}.
\begin{theorem}
\label{T2}
Let the initial data $(u_0,v_0,\theta_0)$ belongs to $H^1(\bR^2)$ with $\dv u_0=0$. Then, the Cauchy problem \eqref{E11} has a unique global solution $(u, v, \te) \in \mathcal{C} ([0,T]; H^1 (\bR^2))$ for all $T>0$. 
\end{theorem}
\begin{proof}
The cited references \cite{MR3479523,MR4969841} establish global existence and uniqueness, together with the following regularity for every $T>0$:
\begin{align*}
(u, v) 
&\in \cC ([0, T];  H^1 (\bR^2)) \cap L^2 (0, T; H^2 (\bR^2)), \\
\theta 
&\in \cC ([0, T];  L^2 (\bR^2)) \cap L^\infty (0, T; H^1 (\bR^2)), \\
(\p_t u, \p_t v, \p_t \theta)
&\in L^2 (0, T; L^2 (\bR^2)), \\
\nabla u& \in L^1 (0, T; L^\infty (\bR^2)).
\end{align*}
It remains to establish the time continuity of $\te$ in $H^1(\bR^2)$.
We provide the proof in Appendix~\ref{app:time_continuity}.
\end{proof}
Theorem~\ref{T2} does not directly yield the preservation of smallness asserted in Theorem~\ref{Thm1}, since it does not provide time-independent bound of
$\|\nabla u\|_{L^1(0,T;L^\infty)}$
(see the estimate in \cite[Corollary~1]{MR3479523}).
This section aims to prove Theorem~\ref{Thm1} in the case $s=1$.
\begin{proposition}
\label{prop_stab_H1}
There exists a sufficiently small constant $\ep_1>0$
with the following property: 
let $(u_0,v_0,\te_0)\in H^1(\bR^2)$ with $\dv u_0=0$.
Suppose that, for some $0<\ep<\ep_1$,
\[
\|(u_0,v_0,\te_0)\|_{H^1 (\bR^2)} < \ep.
\]
\begin{enumerate}[label=(\roman*)]
    \item 
Then the solution obtained in Theorem~\ref{T2} satisfies
the following global-in-time estimate:
\begin{equation}
\label{E12}
\begin{split}
\sup_{t \in [0, \infty)}
\|(u ,v ,\te )(t)\|_{H^1 (\bR^2)}
+ &\left( \int_0^{\infty} 
\|(\nabla u ,\nabla v )(\tau)\|_{H^1}^2
+ \|\nabla \te (\tau)\|_{L^2}^2 \,\dtau \right)^{\frac{1}{2} }
\le 8 \ep.
\end{split}
\end{equation}
\item 
Furthermore, there exists a non-increasing function $\cM_1:[0, \infty)\to[0, \infty)$ 
such that
\begin{equation}
\label{E14}
(1+t)\left(
    \|(\nabla u,\nabla v,\nabla\te)(t)\|_{L^2}^2
    +\int_t^\infty
        \|(\Delta u,\Delta v)(\tau)\|_{L^2}^2
    \,\ud\tau
\right)
\leq 
\cM_1(t)^2,
\qquad \text{for all } t \ge 0,
\end{equation}
where the function $\cM_1$ satisfies
\begin{equation*}
    \sup_{t\in [0, \infty)}\cM_1(t) \leq C\ep,\qquad\textrm{ and }\qquad \lim_{t\to\infty}\cM_1(t)=0. 
\end{equation*}
\end{enumerate}
\end{proposition}

First, we establish an $H^1$ estimate for stability.
For convenience, we state the following lemma.

\begin{lemma}\label{Lem:L2}
Let $(u, v, \te)$ be the strong solution of \eqref{E11} constructed in Theorem~\ref{T2} on $[0, T]$. Then, for all $t\in [0,T]$, we have
\begin{equation*}
\frac{1}{2} 
\|(u, v, \te)(t)\|_{L^2}^2
+ \int_0^t
\norm{(\nabla u, \nabla v)(\tau)}_{L^2}^2
\ud \tau=
\frac{1}{2} 
\|(u_0,v_0,\te_0)\|_{L^2}^2.
\end{equation*}
\end{lemma}
\begin{proof}
We multiply $\eqref{E11}_1$--$\eqref{E11}_3$ by $u,v$ and
$\theta$, respectively. We then integrate over $\mathbb{R}^2$ to obtain
\begin{equation*}
\begin{split}
&\frac{1}{2} \frac{\ud}{\dt}
\|(u, v, \te)\|_{L^2}^2
+ \norm{(\nabla u,\nabla v)}_{L^2}^2
\\
&\quad =- \int_{\bR^2} (u \cdot \nabla) u \cdot u
- \int_{\bR^2}  \dv (v \otimes  v) \cdot u
- \int_{\bR^2} (v \cdot \nabla ) u \cdot v. 
\end{split}
\end{equation*}
By using integration by parts with the divergence-free condition of $u$, we directly obtain
\begin{equation*}
\int_{\bR^2} (u \cdot \nabla) u \cdot u = 0
\qand
\int_{\bR^2} (v \otimes  v) : \nabla u
- \int_{\bR^2} (v \cdot \nabla ) u \cdot v
= 0.
\end{equation*}
\end{proof}

\begin{lemma}
\label{L31}
Let $(u, v, \te)$ be the strong solution of \eqref{E11} constructed in Theorem~\ref{T2} on $[0, T]$.
For each time $t>0$, define $A_1(t)$ and $B_1(t)$ by
\[
A_1 (t) :=  
\|(u, v, \te) (t)\|_{H^1}
\]
and
\[
B_1 (t) := \left( \|(\nabla u,\nabla v) (t)\|_{H^1}^2 
+ \|\nabla \te (t)\|_{L^2}^2 \right)^{\frac{1}{2} }. 
\]
Then, provided that $\|(u_0,v_0,\te_0)\|_{L^2}$ is sufficiently small, we have
\begin{equation}
\label{E34}
A_1 (t)^2 + \frac{1}{2}\int_0^t B_1(\tau)^2 d \tau
\le 3A_1(0)^2
+ 
4\int_0^t 
\norm{\nabla u (\tau)}_{L^\infty}
\norm{\nabla \theta (\tau)}_{L^2}^2 d \tau. 
\end{equation}
\end{lemma}
\begin{proof}
We multiply $\eqref{E11}_1$--$\eqref{E11}_3$ by $(- \Delta) u$, $(- \Delta) v$, and
$(- \Delta) \theta$, respectively, and integrate over $\mathbb{R}^2$ to deduce that
\begin{equation}
\label{E31}
\begin{split}
&\frac{1}{2} \frac{\ud}{\dt}
\left( \norm{\nabla u}_{L^2}^2
+ \norm{\nabla v}_{L^2}^2
+ \|\nabla \te\|_{L^2}^2 
\right)
+ \norm{\Delta u}_{L^2}^2 
+ \|\Delta v\|_{L^2}^2  \\
&= - \int_{\bR^2} \nabla ((u \cdot \nabla) u) : \nabla u
- \int_{\bR^2} \nabla \dv (v \otimes  v) : \nabla u 
- \int_{\bR^2} \nabla ((u \cdot \nabla) v) : \nabla v
- \int_{\bR^2} \nabla ((v \cdot \nabla) u) : \nabla v \\
&\quad - \int_{\bR^2} \nabla (u \cdot \nabla \te) \cdot \nabla \te.
\end{split}
\end{equation}
A main difficulty in the stability analysis of \eqref{E11} is that
the equation for $\theta$ lacks a dissipation term.
To overcome this issue, we introduce the following cross term
involving $v$ and $\theta$.

We multiply $\nabla \theta$ on $\eqref{E11}_2$ and $- \dv v$ on $\eqref{E11}_3$. Then, we integrate over $\mathbb{R}^2$ to obtain 
\begin{align*}
\int_{\bR^2} v_t \cdot \nabla \theta
+ \int_{\bR^2} (u \cdot \nabla) v \cdot \nabla \theta
+ \int_{\bR^2} (v \cdot \nabla) u \cdot \nabla \theta
- \int_{\bR^2} \Delta v \cdot \nabla \theta
&= \int_{\bR^2} | \nabla \theta|^2 \\
\int_{\bR^2} v \cdot \nabla \theta_t 
- \int_{\bR^2} (u \cdot \nabla \theta) \dv v 
&= - \int_{\bR^2} | \dv v|^2.
\end{align*}
After summing these two equations, we use integration by parts to show that
\begin{align*}
\int_{\bR^2} |\nabla \theta|^2
&+ \int_{\bR^2} \Delta v \cdot \nabla \theta
- \int_{\bR^2} |\dv v|^2 
- \frac{\ud}{\dt} \int_{\bR^2} v \cdot \nabla \theta \\
&= \int_{\bR^2} (u \cdot \nabla) v \cdot \nabla \theta
+ \int_{\bR^2} (v \cdot \nabla) u \cdot \nabla \theta
- \int_{\bR^2} ( u \cdot \nabla \theta) \cdot \dv v= 2 \int_{\bR^2}  (v \cdot \nabla) u \cdot \nabla \theta.
\end{align*}
We then use Young's inequality to see 
\begin{equation}
\label{E32}
\frac{1}{2} \|\nabla \te\|_{L^2}^2
- \|\nabla v\|_{H^1}^2
- \frac{\ud}{\dt} \int_{\bR^2} v \cdot \nabla \theta \le  2 \int_{\bR^2}  (v \cdot \nabla) u \cdot \nabla \theta.
\end{equation}
Finally, we multiply \eqref{E32} by a sufficiently small constant, for instance $1/2$,
and add it to \eqref{E31}. We also combine Lemma~\ref{Lem:L2} to obtain
\begin{equation}
\label{E33}
\begin{split}
&\frac{1}{2} \frac{\ud}{\dt}
\left( \| \nabla u\|_{L^2}^2
+ \|\nabla v\|_{L^2}^2
+ \|\nabla \te\|_{L^2}^2 
+ \|u\|_{L^2}^2
+ \|v\|_{L^2}^2
+\| \te \|_{L^2}^2
- \int_{\bR^2} v \cdot \nabla \theta
\right) \\
&\quad + \| \Delta u\|_{L^2}^2 
+ \frac{1}{2} \|\Delta v\|_{L^2}^2
+ \| \nabla u\|_{L^2}^2
+ \frac{1}{2} \|\nabla v \|_{L^2}^2
+ \frac{1}{4} \|\nabla \te\|_{L^2}^2 \\
&
\le -\int_{\bR^2} \nabla ((u \cdot \nabla) u) : \nabla u
- \int_{\bR^2} \nabla \divg (v \otimes  v) : \nabla u 
- \int_{\bR^2} \nabla ((u \cdot \nabla) v) : \nabla v
- \int_{\bR^2} \nabla ((v \cdot \nabla) u) : \nabla v\\    
&\quad - \int_{\bR^2} \nabla (u \cdot \nabla \te) \cdot \nabla \te
+ \int_{\bR^2}  (v \cdot \nabla) u \cdot \nabla \theta\\
&=: \sum_{k = 1}^6\cI_k.
\end{split}
\end{equation}

We now estimate the terms $\mathcal I_k$ for $k = 1, \dots, 6$, on the right-hand side of \eqref{E33}.
We simply refer to \cite[Lemma 2.1]{BCK2025} to reveal that 
\begin{equation*}
    \cI_1 = \int_{\bR^2} \nabla ((u \cdot \nabla) u) \cdot \nabla u = \sum_{i,j,k = 1,2} \int_{\bR^2}\partial_i u_j \partial_ju_k \partial_i u_k=0.
\end{equation*}

Applying the Gagliardo–Nirenberg interpolation inequality leads to
\begin{equation*}
\|w\|_{L^4 (\bR^2)} \lesssim \norm{ w}_{L^2 (\bR^2)}^\frac{3}{4} 
\norm{\Delta w}_{L^2 (\bR^2)}^\frac{1}{4} \quad \text{and} \quad 
\|\nabla w\|_{L^{4} (\bR^2)} \lesssim \norm{ w}_{L^2 (\bR^2)}^\frac{1}{4} 
\norm{\Delta w}_{L^2 (\bR^2)}^\frac{3}{4}. 
\end{equation*}

We then use these together with a standard commutator estimate, under the divergence-free condition of $u$, to obtain
\begin{equation}\label{eq:smallness1}
\begin{aligned}
|\cI_2|
&= \lt|\,\int_{\bR^2} \divg (v \otimes  v) \cdot \Delta u \,\rt|
\lesssim \norm{v}_{L^4}
\norm{\nabla v}_{L^4}
\norm{\Delta u}_{L^2}
\lesssim
\norm{v}_{L^2}
\norm{\Delta v}_{L^2} 
\norm{\Delta u}_{L^2}, \\
|\cI_3|
&=\lt|\,\int_{\bR^2} (u\cdot \nabla v) \cdot \Delta v \,\rt|
\lesssim \norm{u}_{L^4}
\norm{\nabla v}_{L^4}
\norm{\Delta v}_{L^2}
\lesssim
\norm{u}_{L^2}^\frac{3}{4}
\norm{v}_{L^2}^\frac{1}{4}
\norm{\Delta v}_{L^2}^\frac{7}{4} 
\norm{\Delta u}_{L^2}^\frac{1}{4}, \\
|\cI_4|
&=\lt|\,\int_{\bR^2} (v\cdot \nabla u) \cdot \Delta v \,\rt|
\lesssim \norm{v}_{L^4}
\norm{\nabla u}_{L^4}
\norm{\Delta v}_{L^2}
\lesssim
\norm{v}_{L^2}^\frac{3}{4}
\norm{u}_{L^2}^\frac{1}{4}
\norm{\Delta u}_{L^2}^\frac{3}{4} 
\norm{\Delta v}_{L^2}^\frac{5}{4}, \\
|\cI_5| 
&\le \norm{\nabla u}_{L^\infty}
\|\nabla \te\|_{L^2}^2, \\
|\cI_6|
&\le 
\norm{v}_{L^ 4} 
\norm{\nabla u}_{L^4}
\norm{\nabla \theta}_{L^2}
\lesssim
\norm{v}_{L^2}^\frac{3}{4}
\norm{u}_{L^2}^\frac{1}{4}
\norm{\Delta v}_{L^2}^\frac{1}{4} 
\norm{\Delta u}_{L^2}^\frac{3}{4}
\norm{\nabla \theta}_{L^2}.
\end{aligned}
\end{equation}
It follows from Lemma~\ref{Lem:L2} that
\begin{equation}\label{eq:smallness2}
    \begin{split}
|\cI_2|+|\cI_3|+|\cI_4|
&\le C (\norm{u}_{L^2}
+ \norm{v}_{L^2})
(\norm{\Delta u}_{L^2}^2
+ \norm{\Delta v}_{L^2}^2)
\leq C\|(u_0,v_0,\te_0)\|_{L^2}B_1^2(t)
\end{split}
\end{equation}
and
\begin{equation*}
    |\cI_6| \leq C\|(u_0,v_0,\te_0)\|_{L^2}B_1^2(t).
\end{equation*}
Collecting the estimates for $\cI_k$, we obtain
\[
\sum_{k=1}^6 |\cI_k|
\leq
\widetilde{C}\|(u_0,v_0,\te_0)\|_{L^2} B_1 (t)^2
+ \norm{\nabla u}_{L^\infty}
\|\nabla \te\|_{L^2}^2.
\]
Choosing the initial $L^2$-norm $\|(u_0,v_0,\te_0)\|_{L^2}$ sufficiently small to ensure
$\widetilde{C}\|(u_0,v_0,\te_0)\|_{L^2}\leq \frac{1}{8}$ and absorbing the first term into the left-hand side, we have
\begin{equation}\label{eq0526wed}
\frac{\ud}{\dt}
\tilde{A}_1^2 
+ \frac{1}{4} B_1^2
\leq 2\norm{\nabla u}_{L^\infty}
\|\nabla \te\|_{L^2}^2,
\end{equation}
where
\begin{equation*}
    \tilde{A}_1(t) :=
    \left(
\|u(t)\|_{H^1}^2
+ \|v(t)\|_{H^1}^2
+\| \te (t)\|_{H^1}^2
- \int_{\bR^2} v(t) \cdot \nabla \theta(t)
\right)^{\frac{1}{2}}.
\end{equation*}
It is easy to verify that
\[
\frac{1}{2}
A_1(t)^2
\le  
\tilde{A}_1(t)^2
\le \frac{3}{2} 
A_1(t)^2.
\]
Using this equivalence, we integrate \eqref{eq0526wed} over $[0,t]$ to obtain
\[
\frac{1}{2} A_1 (t)^2 
+ \frac{1}{4}\int_0^t B_1 (\tau)^2 d \tau
\le \frac{3}{2} A_1 (0)^2  
+ 2\int_0^t \norm{\nabla u (\tau)}_{L^\infty}
\norm{\nabla \theta (\tau)}_{L^2}^2 d \tau. 
\]
The proof is complete.
\end{proof}

\medskip
\noindent
It remains to estimate the term
\begin{equation*}
\int_0^t \|\nabla u (\tau)\|_{L^\infty}
\|\nabla \theta (\tau)\|_{L^2}^2 d \tau.
\end{equation*}
To this end, we first introduce the weighted interpolation inequality. The proof is a slight modification of the argument in \cite[Lemma 2.1]{JKL2022}, and thus we omit the proof.

\begin{lemma}
\label{L23}
Let $1 \le q < 2$ and let $a,b > 0$. Then there exists a constant $C = C (a,b,q) > 0$ such that
\begin{equation*}
\norm{f}_{L^q (\bR^2)}
\le C \norm{|\cdot|^{\frac{2}{q} - 1 + a} f}_{L^2 (\bR^2)}^\frac{b}{a+b}
\norm{|\cdot|^{\frac{2}{q} - 1 -b} f}_{L^2 (\bR^2)}^\frac{a}{a+b},
\end{equation*}
for all $f \in \mathscr{S} (\bR^2)$.
In particular, if $1 < q < 2$, choosing
$a=2- \frac{2}{q}$ and $b=\frac{2}{q}-1$ yields
\begin{equation*}
    \norm{\sF u}_{L^q(\bR^2)}
\lesssim_{q}
\norm{\nabla u}_{L^2 (\bR^2)}^\frac{2-q}{q}
\norm{u }_{L^2 (\bR^2)}^\frac{2q - 2}{q}.
\end{equation*}
\end{lemma}

Observing that
\[
\|\nabla u (t)\|_{L^\infty}
\le C \int_{\bR^2} |\xi| |\sF u (t, \xi)| d \xi,
\]
we estimate $\|\,|\xi|\sF u (t, \xi)\|_{L^1_\xi}$ via the Duhamel formula.

\begin{lemma}
\label{Prop1}
Let $1<p<2$, and let $r$ and $t$ satisfy $0 \le r < t \le T < \infty$.
Let $(u, v, \te)$ be a strong solution of \eqref{E11} on
$[0, T] \times \bR^2$, with
$u, v \in C ([0, T]; H^1(\bR^2)) \cap L^2 (0, T; H^2(\bR^2))$.
Then,
\begin{equation}
\label{eq:modi_1}
\begin{split}
\left( \int_r^t\norm{\nabla u (\tau)}_{L^\infty}^p d \tau \right)^{\frac{1}{p}}
&\lesssim_{p}
\|u(r)\|_{L^2}^{\frac{2}{p} - 1} \|\nabla u(r)\|_{L^2}^{2 - \frac{2}{p}}
+ \left(
\sup_{\tau \in [r, t]} 
\norm{u (\tau)}_{L^2}
\right)^{1 - \frac{1}{p}}
\norm{\nabla u}_{L^2 ([r, t]; L^2)}^{\frac{3}{p} - 1}
\norm{\Delta u}_{L^2 ([r, t]; L^2)}^{2 - \frac{2}{p}} \\
&\qquad + \left(
\sup_{\tau \in [r, t]} 
\norm{v (\tau)}_{L^2}
\right)^{1 - \frac{1}{p}}
\norm{\nabla v}_{L^2 ([r, t]; L^2)}^{\frac{3}{p} - 1}
\norm{\Delta v}_{L^2 ([r, t]; L^2)}^{2 - \frac{2}{p}}.
\end{split}
\end{equation}
Here, the implicit constant depends only on $p$, and is independent of $r,t,T>0$ and the solution $(u, v, \te)$.
\end{lemma}

\begin{proof}
Without loss of generality, we assume that $r=0$ and $t=T$.
Applying the Fourier transform $\sF$ to $\eqref{E11}_1$, we obtain
\[
\sF u_t  
+ |\xi|^2 \sF u
= - \mathbb P (\xi) \sF ((u \cdot \nabla) u) 
- \mathbb P (\xi) \sF \dv (v \otimes v),
\]
where $\mathbb P (\xi)
= I - \xi \otimes \xi / |\xi|^2$, ($\xi \neq 0$) denotes the Leray projection.
The Duhamel formula then yields
\[
\sF u (t, \xi)
= e^{- |\xi|^2 t} \sF u_0 (\xi)
- \int_0^t e^{- |\xi|^2 (t - \tau)} \mathbb P (\xi) \sF ((u \cdot \nabla) u) (\tau, \xi) d \tau
- \int_0^t e^{- |\xi|^2 (t - \tau)} \mathbb P (\xi) \sF \dv (v \otimes v) (\tau, \xi) d \tau.
\]
Using that $|\mathbb P (\xi)| \le 1$, we have
\begin{equation}
\label{E211}
\int_{\bR^2} |\xi| |\sF u (t, \xi)| d \xi
\le \cJ_0 (t) 
+ \cJ_1 (t) 
+ \cJ_2 (t),
\end{equation}
where
\begin{align*}
\cJ_0 (t) 
&:= \int_{\bR^2} |\xi| e^{- |\xi|^2 t} |\sF u_0 (\xi)| d \xi, \\
\cJ_1 (t)
&:= \int_{\bR^2} |\xi|\int_0^t e^{- |\xi|^2 (t - \tau)} 
|\sF ((u \cdot \nabla) u) (\tau, \xi)| d \tau d \xi, \\
\cJ_2 (t)
&:= \int_{\bR^2} |\xi|\int_0^t e^{- |\xi|^2 (t - \tau)} |\sF \dv (v \otimes v) (\tau, \xi)| d \tau d \xi.
\end{align*}

We first estimate the first integrand $\cJ_0$ on the right-hand side of \eqref{E211}.
By Minkowski’s integral inequality, we have
\begin{equation*}
\begin{split}
\|\cJ_0\|_{L^p_t(0, T)}
&=\Norm{\int_{\bR^2} |\xi| e^{- |\xi|^2 t} |\sF u_0 (\xi)| d \xi}_{L_t^p (0, T)}
\le \int_{\bR^2} \left( \int_0^T |\xi|^p e^{- tp |\xi|^2} |\sF u_0 (\xi)|^p dt \right)^\frac{1}{p} d \xi \\
&\le \int_{\bR^2} |\xi| |\sF u_0 (\xi)| \left( \int_0^\infty e^{- tp |\xi|^2} dt \right)^\frac{1}{p} d \xi
= C_p \int_{\bR^2} |\xi|^{1 - \frac{2}{p}} |\sF u_0 (\xi)| d \xi.
\end{split}
\end{equation*}

Since $p \in (1, 2)$, for $R>0$, we have
\begin{align*}
\int_{\bR^2} |\xi|^{1 - \frac{2}{p}} |\sF u_0 (\xi)| d \xi 
&= \int_{|\xi| < R} |\xi|^{1 - \frac{2}{p}} |\sF u_0 (\xi)| d \xi
+ \int_{|\xi| \ge R} |\xi|^{1 - \frac{2}{p}} |\sF u_0 (\xi)| d \xi \\
&\le \left( \int_{|\xi| < R} |\xi|^{2 \left(1 - \frac{2}{p} \right)} d \xi \right)^\frac{1}{2} 
\left( \int_{|\xi| < R} |\sF u_0 (\xi)|^2 d \xi \right)^\frac{1}{2} 
+ \left( \int_{|\xi| \ge R} |\xi|^{- \frac{4}{p} } d \xi \right)^\frac{1}{2} 
\left( \int_{|\xi| \ge R} |\xi|^2|\sF u_0 (\xi)|^2 d \xi \right)^\frac{1}{2} \\
&\le C  R^{2 - \frac{2}{p}} \|u_0\|_{L^2}
  +CR^{1- \frac{2}{p}} \|\nabla u_0\|_{L^2},
\end{align*}
where the constant depends only on $p$.

If $u_0 = 0$, then $\cJ_0 = 0$. Otherwise, by taking $R= \|\nabla u_0\|_{L^2} /\|u_0\|_{L^2}$, we have
\begin{equation}
\label{E49}
\|\cJ_0\|_{L^p_t(0, T)}= \Norm{\int_{\bR^2} |\xi| e^{- |\xi|^2 t} |\sF u_0 (\xi)| d \xi}_{L_t^p (0, T)}
\le C_p \|u_0\|_{L^2}^{\frac{2}{p} - 1} \|\nabla u_0\|_{L^2}^{2 - \frac{2}{p}} .
\end{equation}

Now we estimate the second integrand $\cJ_1$ on the right-hand side of \eqref{E211}.
One applies Tonelli's theorem and Young's convolution inequality to obtain
\begin{equation}
\label{E213}
\begin{split}
\cJ_1 (t) 
= \int_{\bR^2} |\xi| \left( \int_0^t e^{- |\xi|^2 (t - \tau)} |\sF ((u \cdot \nabla) u) (\tau, \xi)| d \tau \right) d \xi 
&\le C \int_0^t \int_{\bR^2} 
|\xi| e^{- |\xi|^2 (t - \tau)} 
\left( |\sF u| \ast |\sF \nabla u| \right) (\tau, \xi) d \xi d \tau \\
&\le C \int_0^t \norm{|\cdot| e^{- (t - \tau) |\cdot|^2}}_{L_\xi^r}
\norm{\sF u (\tau)}_{L_\xi^{q_1}}
\norm{\sF \nabla u (\tau)}_{L_\xi^{q_2}} d \tau,
\end{split}
\end{equation}
where the exponents $(r,q_1,q_2)$ are defined by
\begin{equation}
\label{Exponent}
    r:=2p \in (2,4), \quad q_1:=\frac{2p}{p+1} \in (1, 2), \quad q_2:=\frac{2p}{3p-2} \in (1, 2).
\end{equation}
Note that the exponents satisfy
 $1+ \frac{1}{r/(r-1)}=\frac{1}{q_1} + \frac{1}{q_2}$.
Using the change of variables, one finds that
\begin{equation}
\label{E215}
\norm{|\cdot| e^{- (t - \tau) |\cdot|^2}}_{L_\xi^r}^r
= \int_{\bR^2} |\xi|^r e^{- r (t - \tau) |\xi|^2} d \xi
= \int_{\bR^2} \left|\frac{\eta}{\sqrt{r (t - \tau)}} \right|^r e^{- r (t - \tau) \left|\frac{\eta}{\sqrt{r (t - \tau)}} \right|^2} \frac{d \eta}{r (t - \tau)}
\lesssim (t - \tau)^{- \frac{r + 2}{2}}.
\end{equation}
Since $q_1, q_2 \in (1, 2)$, we can apply Lemma \ref{L23} to obtain that 
\begin{equation}
\label{E216}
\norm{\sF u (\tau)}_{L_\xi^{q_1}}
\lesssim_{p}
\norm{\nabla u (\tau)}_{L^2}^\frac{2-q_1}{q_1}
\norm{u (\tau)}_{L^2}^\frac{2q_1 - 2}{q_1}
\quad \text{and} \quad
\norm{\sF \nabla u (\tau)}_{L_\xi^{q_2}}
\lesssim_{p} 
\norm{\Delta u (\tau)}_{L^2}^\frac{2-q_2}{q_2}
\norm{\nabla u (\tau)}_{L^2}^\frac{2q_2 - 2}{q_2}.
\end{equation}
Combining \eqref{E213} with \eqref{E215}-\eqref{E216}, it follows that

\begin{equation}
\label{E45}
\cJ_1 (t) 
\lesssim_{p}
\int_0^t 
(t - \tau)^{- \frac{r + 2}{2r}}
\norm{\Delta u (\tau)}_{L^2}^{\frac{2-q_2}{q_2}}
\norm{\nabla u (\tau)}_{L^2}^{1 +\frac{2}{q_1}- \frac{2}{q_2}}
\norm{u (\tau)}_{L^2}^{\frac{2q_1-2}{q_1}}
d \tau.
\end{equation}

For the exponents defined in \eqref{Exponent}, the
Hardy--Littlewood--Sobolev inequality gives
\begin{equation}
\label{E410}
\Norm{\int_0^t (t - \tau)^{- \frac{r + 2}{2 r}} f (\tau) d \tau}_{L_t^p (0, T)}
\lesssim \norm{f}_{L_t^{q_1} (0, T)}.
\end{equation}
Indeed, the exponents satisfy $1 < q_1 < p < 2$ and $\frac{1}{p} = \frac{1}{q_1} - \frac{r - 2}{2 r}$ with $0<\frac{r-2}{2r}<1$.
The implicit constant is independent of $T$.

Consequently, we directly apply \eqref{E410} to \eqref{E45} and deduce 
\begin{equation}
\begin{split}
\label{E411}
\|\cJ_1\|_{L^p_t(0, T)} 
&\lesssim \left\|
\norm{\Delta u (\tau)}_{L^2}^{\frac{2-q_2}{q_2}}
\norm{\nabla u (\tau)}_{L^2}^{1 +\frac{2}{q_1}- \frac{2}{q_2}}
\norm{u (\tau)}_{L^2}^{\frac{2q_1-2}{q_1}}
\right\|_{L^{q_1}_t(0, T)}
\\
&\le \sup_{t \in [0, T]} \norm{u (t)}_{L^2}^{\frac{2q_1-2}{q_1}}
\left\|   
\| \nabla u (\tau) \|_{L^2}^{1 +\frac{2}{q_1}- \frac{2}{q_2}}
\right\|_{L^{2p/(3-p)}_t(0, T)}
\left\|   
\| \Delta u (\tau) \|_{L^2}^{\frac{2-q_2}{q_2}}
\right\|_{L^{p/(p-1)}_t(0, T)} \\
&\le \left(
\sup_{t \in [0, T]} 
\norm{u(t)}_{L^2}
\right)^{1 - \frac{1}{p}}
\norm{\nabla u}_{L^2_t ([0, T]; L^2)}^{\frac{3}{p} - 1}
\norm{\Delta u}_{L^2_t ([0, T]; L^2)}^{2 - \frac{2}{p}}.
\end{split}
\end{equation}
A similar estimate can be obtained for $\|\cJ_2\|_{L^p_t(0, T)}$ by repeating the above argument:
\begin{equation}
\label{E412}
\|\cJ_2\|_{L^p_t(0, T)} 
\lesssim 
\left(\sup_{t \in [0, T]} \norm{v(t)}_{L^2} \right)^{1 - \frac{1}{p}}
\norm{\nabla v}_{L^2_t ([0, T]; L^2)}^{\frac{3}{p} - 1}
\norm{\Delta v}_{L^2_t ([0, T]; L^2)}^{2 - \frac{2}{p}}.
\end{equation}
Finally, combining 
\eqref{E49}, \eqref{E411}, and \eqref{E412} together with \eqref{E211} leads to the desired result.
\end{proof}

\medskip
Now, we are ready to prove Proposition~\ref{prop_stab_H1}.

\begin{proof}[Proof of the Proposition~\ref{prop_stab_H1}]
We divide the proof into two parts: $(i)$ \textbf{Stability.} $(ii)$ \textbf{Temporal decay.}

\medskip
\noindent
\emph{(i) \textbf{Stability.}}

We first establish \eqref{E12}. 
Define a functional
\[
E (t):=  \sup_{\tau \in [0, t]} 
\| (u, v, \te) (\tau) \|_{H^1}^2
+ \frac{1}{2}\int_0^t \left( \|(\nabla u, \nabla v) (\tau)\|_{H^1}^2 
+  \|\nabla \te (\tau)\|_{L^2}^2 \right) d \tau.
\]
Then, by Lemma~\ref{L31}, we obtain that for $t > 0$,
\begin{equation}
\label{E50}
E (t)
\le 6 E (0)
+ 8 \int_0^t \norm{\nabla u (\tau)}_{L^\infty}
\norm{\nabla \theta (\tau)}_{L^2}^2 d \tau.
\end{equation}
We fix $p \in (1, 2)$. 
By H\"older's inequality,
\begin{equation}
\label{E51}
\int_0^t \norm{\nabla u (\tau)}_{L^\infty}
\norm{\nabla \theta (\tau)}_{L^2}^2 d \tau
\le \sup_{\tau \in [0, t]} 
\norm{\nabla \theta (\tau)}_{L^2}^{\frac{2}{p}}
\left( \int_0^t \norm{\nabla u (\tau)}_{L^\infty}^p d \tau \right)^{\frac{1}{p}}
\left( \int_0^t \norm{\nabla \te (\tau)}_{L^2}^2 d \tau \right)^{1 - \frac{1}{p}}.
\end{equation}
From the definition of $E (t)$, we immediately have
\begin{equation}
\label{E52}
\left( \sup_{\tau \in [0, t]} 
\norm{\nabla \theta (\tau)}_{L^2} \right)^{\frac{2}{p}}
\norm{\nabla \theta}_{L^2 ([0, t]; L^2)}^{2 - \frac{2}{p}}
\le E(t)^{1/p}(2E(t))^{1-1/p}
=2^{1-1/p} E(t).
\end{equation}
On the other hand, Lemma~\ref{Prop1}, with $r=0$, implies that
\begin{equation}
\label{E53}
\left( \int_0^t\norm{\nabla u (\tau)}_{L^\infty}^p d \tau \right)^{\frac{1}{p}}
\lesssim_p E (0)^\frac{1}{2} 
+ \left( E (t)^\frac{1}{2} \right)^{1 - \frac{1}{p}}
\left( E (t)^\frac{1}{2} \right)^{1 + \frac{1}{p}} 
=
E (0)^\frac{1}{2} 
+ E (t).
\end{equation}
Combining \eqref{E50} and \eqref{E51}-\eqref{E53} with Young's inequality, we obtain
\begin{equation}
\label{E54}
E (t)
\le 7 E (0)
+ C E (t)^2.
\end{equation}

By Theorem \ref{T2}, the functional $E (t)$ is continuous.
By a standard bootstrapping argument, we conclude from \eqref{E54} that there exists a constant $\ep_1 > 0$ such that if $0 < \ep < \ep_1$ and $E (0) < \ep$, then 
\[
E (t) \le 8 E (0), \quad \text{for all } t \in [0, \infty).
\]
This proves \eqref{E12}.

\medskip
\noindent
\emph{(ii) \textbf{Temporal decay.}}

We next prove the temporal decay estimate \eqref{E14} using the stability estimate \eqref{E12}.
For sufficiently small $H^1$-norm $\|(u_0,v_0,\te_0)\|_{H^1}>0$,
the estimate \eqref{E31} with
\eqref{eq:smallness1}--\eqref{eq:smallness2} yields
\begin{equation}
\label{E63}
\begin{split}
&\frac{\ud}{\dt}
\norm{(\nabla u, \nabla v, \nabla \te)}_{L^2}^2
+ \norm{(\Delta u, \Delta v)}_{L^2}^2
\le 2 \norm{\nabla u}_{L^\infty}
\norm{\nabla \theta}_{L^2}^2.
\end{split}
\end{equation}

Fix $1<p<2$.
Integrating \eqref{E63} over $[r,\rho]$ and $[r,t]$, and then taking the supremum over $\rho\in[r,t]$ and and adding the resulting inequalities, we find that
\begin{equation}
\label{K11}
\begin{split}
&\sup_{\rho \in [r, t]}
\norm{
(\nabla u, \nabla v, \nabla \te) (\rho)}_{L^2}^2
+ \int_r^t
\norm{(\Delta u, \Delta v) (\tau)}_{L^2}^2
\, \dtau \\
&\le 2
\norm{(\nabla u, \nabla v, \nabla \te) (r)}_{L^2}^2
+ 4 \sup_{\rho \in [r,t]} 
\norm{\nabla \theta (\rho)}_{L^2}^{\frac{2}{p}}
\left( \int_r^t \norm{\nabla u (\tau)}_{L^\infty}^p d \tau \right)^{\frac{1}{p}}
\left( \int_r^t \norm{\nabla \te (\tau)}_{L^2}^2 d \tau \right)^{1 - \frac{1}{p}}
\end{split}
\end{equation}
for every $0\le r<t$.
We used H\"older's inequality.

For each fixed $r\geq0$, we define a functional
\[
G_r(t):=
\sup_{\tau\in[r,t]}
    \|(\nabla u
    , \nabla v, \nabla \te)(\tau)\|_{L^2}^2
+
\int_r^t
    \|(\Delta u, \Delta v)(\tau)\|_{L^2}^2\,\ud\tau,
\qquad t\geq r.
\]
For every $0 \le r < t$, Lemma~\ref{Lem:L2} and \eqref{E12} give
\begin{equation}\label{K17}
    \sup_{\tau \in [r, t]} \norm{(u,v) (\tau)}_{L^2}
\le \|(u_0,v_0,\te_0)\|_{L^2}
\le \ep
\quad \text{and} \quad
\int_r^t 
\norm{(\nabla u,
\nabla v,
\nabla \te)(\tau)}_{L^2}^2  \,\dtau
\le 64 \ep^2.
\end{equation}
Applying Lemma~\ref{Prop1} on $[r,t]$, we therefore have
\begin{equation}
\label{K12}
\begin{aligned}
\left(
    \int_r^t\|\nabla u(\tau)\|_{L^\infty}^p\,\ud\tau
\right)^{1/p}
&\lesssim_p
\ep^{2/p-1}
\|\nabla u(r)\|_{L^2}^{2-2/p}
+
\ep^{2/p}
\left(
    \int_r^t
        \|(\Delta u, \Delta v)(\tau)\|_{L^2}^2
        \,\ud\tau
\right)^{1-1/p} \\
&\lesssim_p
    \ep^{2/p-1}G_r(r)^{1-1/p}
    +\ep^{2/p}G_r(t)^{1-1/p}.
\end{aligned}
\end{equation}

Then, \eqref{K11}-\eqref{K12} implies that, for $t \ge r$, 
\begin{align*}
G_r(t)
&\leq
2G_r(r)
+C G_r(t)^{1/p}\ep^{2-2/p}
\left(
    \int_r^t\|\nabla u(\tau)\|_{L^\infty}^p\,\ud\tau
\right)^{1/p}
\\
&\leq
2G_r(r)
+C G_r(t)^{1/p}\ep^{2-2/p}
\left(
    \ep^{2/p-1}G_r(r)^{1-1/p}
    +\ep^{2/p}G_r(t)^{1-1/p}
\right)
\\
&\leq
4G_r(r)+C_1(\ep^p+\ep^2)G_r(t),
\end{align*}
where we used Young's inequality in the last line.

We take $\ep_1$ sufficiently small so that $C_1 (\ep_1^2 + \ep_1^p) \le1/2$ to obtain
\begin{equation}\label{eq1016wed}
G_r(t) 
\le 8G_r(r) = 8 \left(
\norm{\nabla u (r)}_{L^2}^2
+ \norm{\nabla v (r)}_{L^2}^2
+ \norm{\nabla \theta (r)}_{L^2}^2
\right).
\end{equation}
Now, we use a minimum-selection argument. Fix $T \geq1$.
By the time-continuity, the map $t\mapsto \norm{(\nabla u, \nabla v, \nabla \te)(t)}_{L^2}^2 $ attains its minimum on $[T/4, T/2]$.
Choose a minimizer $r^* \in [T/4,T/2]$. Then, it follows that
\begin{equation}
\label{E64}
\begin{split}
\norm{(\nabla u, \nabla v, \nabla \te)(r^*)}_{L^2}^2
&\le \frac{4}{T} \int_{T/4}^{T/2}
\norm{(\nabla u, \nabla v, \nabla \te)(\tau)}_{L^2}^2 \dtau\\
&\le \frac{4}{T} \int_{T/4}^{\infty}
\norm{(\nabla u, \nabla v, \nabla \te)(\tau)}_{L^2}^2  \dtau
\le \frac{256 \ep^2}{T},
\end{split}
\end{equation}
where the last inequality follows from \eqref{E12}.
Taking $r=r^*$ and $t=T$ to \eqref{eq1016wed}, and using
\eqref{E64} together with $r^*\leq T/2$, we conclude that
\begin{equation}
\begin{gathered}\label{E66}
\norm{(\nabla u, \nabla v, \nabla \te)(T)}_{L^2}^2  
+
\int_{T/2}^T
    \|(\Delta u, \Del v)(\tau)\|_{L^2}^2\,\ud\tau
\\
\leq G_{r^*}(T)
\leq8G_{r^*}(r^*)
\leq\frac{\cM_1(T)^2}{T},
\end{gathered}
\end{equation}
where
\[
    \cM_1(T):=
    8\left(
        \int_{T/4}^\infty
       \norm{(\nabla u, \nabla v, \nabla \te)(\tau)}_{L^2}^2\,\ud\tau
    \right)^{1/2}.
\]
By \eqref{E12},
\[
    \sup_{t\in [0, \infty)}\cM_1(t)\leq64\ep
    ,
    \qquad\lim_{T\to\infty}\cM_1(T)=0.
\]
Moreover, by observing that $\cM_1$ is non-increasing,
we apply \eqref{E66} with $T$ replaced by $2^{k+1}T$
and summing over $k \ge 0$ to deduce that
\begin{align*}
\int_T^\infty
\|(\Delta u, \Del v)(\tau)\|_{L^2}^2\,\ud\tau 
&\leq
\sum_{k=0}^\infty
\int_{2^kT}^{2^{k+1}T}
\|(\Delta u, \Del v)(\tau)\|_{L^2}^2\,\ud\tau \\
&\leq
\sum_{k=0}^\infty
\frac{\cM_1(2^{k}T)^2}{2^{k}T}
\le
\frac{\cM_1(T)^2}{T}
\sum_{k=0}^\infty 2^{-k}
=
\frac{2\cM_1(T)^2}{T}.
\end{align*}
Combining this estimate with \eqref{E66}, we conclude that
\[
\begin{aligned}
&\norm{(\nabla u, \nabla v, \nabla \te)(T)}_{L^2}^2
+
\int_T^\infty
    \|(\Delta u, \Del v)(\tau)\|_{L^2}^2\,\ud\tau
\le \frac{3\cM_1(T)^2}{T},
\qquad \textrm{for }T \geq1.
\end{aligned}
\]
Redefining $\cM_1$ by a constant factor, we obtain \eqref{E14}.
The proof is complete.
\end{proof}

\subsection{Time-integral estimates on dyadic intervals}
The $\dot{H}^1$ temporal decay estimates established above yield the following Lemma~\ref{lem:dyadic_int}. The uniform bounds on dyadic-time integrals in Lemma~\ref{lem:dyadic_int} will be used in the proof of the temporal decay part of the main theorem.
As suggested by \eqref{post:1707} and \eqref{post:1708},
these bounds appear difficult to obtain without the
$\dot{H}^1$ decay estimates established above.
\begin{lemma}
\label{lem:dyadic_int}
Under the assumptions of Proposition~\ref{prop_stab_H1},
the solution $u$ satisfies
\[
    \int_T^{2T}
    \left(
        \|\nabla u(t)\|_{L^\infty}
        +\|\Delta u(t)\|_{L^2}
    \right)\,\ud t
    \leq C\ep
    \qquad (T\geq1),
\]
where $C>0$ is an absolute constant.
\end{lemma}

\begin{proof}
Recall from \eqref{E14} and its bound $\sup \cM_1 \leq C\ep$ that
\[
    \|\nabla u(T)\|_{L^2}^2
    +\int_T^{2T}
        \|(\Delta u, \Delta v)(t)\|_{L^2}^2\,\ud t
    \leq \frac{C\ep^2}{T}.
\]

Let $1<p<2$. 
We apply \eqref{K12} on $[T,2T]$:
\[
\left(
    \int_T^{2T}\|\nabla u(\tau)\|_{L^\infty}^p\,\dtau
\right)^{1/p}
\lesssim_p
\ep^{2/p-1}
\|\nabla u(T)\|_{L^2}^{2-2/p}
+
\ep^{2/p}
\left(
    \int_T^{2T}
        \|(\Delta u,\Delta v)(\tau)\|_{L^2}^2
\,\ud\tau
\right)^{1-1/p}.
\]
Combining these, we find that
\begin{equation}\label{post:1707}
\begin{aligned}
    \int_T^{2T}\|\nabla u(t)\|_{L^\infty}\,\ud t
    &\leq
    T^{1-1/p}
    \left(
        \int_T^{2T}\|\nabla u(t)\|_{L^\infty}^p\,\ud t
    \right)^{1/p}
\\
    &\leq
    C_p T^{1-1/p}
    \left(
        \ep^{2/p-1}
        \left(\frac{\ep^2}{T}\right)^{1-1/p}
        +\ep^{2/p}
        \left(\frac{\ep^2}{T}\right)^{1-1/p}
    \right)
\\
    &=
    C_p(\ep+\ep^2).
    \end{aligned}
\end{equation}
Fixing $p=3/2$ and assuming $\ep_1\leq1$, we conclude that
\[
    \int_T^{2T}\|\nabla u(t)\|_{L^\infty}\,\ud t
    \leq C\ep,
\]
where $C$ is an absolute constant independent of $T\geq1$.
Similarly, the Cauchy--Schwarz inequality gives
\begin{equation}\label{post:1708}
    \int_T^{2T}\|\Delta u(t)\|_{L^2}\,\ud t
    \leq
    T^{1/2}
    \left(
        \int_T^{2T}\|\Delta u(t)\|_{L^2}^2\,\ud t
    \right)^{1/2}
    \leq 
    T^{\frac{1}{2}}
    \left(\frac{C\ep^2}{T}\right)^{\frac{1}{2}}
    \leq C \ep.
\end{equation}
This completes the proof.
\end{proof}

\section{Higher regularity: $s>1$}
\label{S3}
We begin with the following local existence and uniqueness result for strong solutions in $H^s$.
\begin{theorem}[Local-in-time existence]
\label{T3}
Let $s>1$, and let the initial data $(u_0,v_0,\te_0)$ belong to
$H^s(\bR^2)$ with $\dv u_0=0$.
Then there exists $T>0$ such that the Cauchy problem \eqref{E11}
admits a unique strong solution satisfying
\[
    (u, v, \te)\in\cC([0,T];H^s(\bR^2)),
    \qquad
    (u,v)\in L^2(0,T;H^{s+1}(\bR^2)).
\]
\end{theorem}
\begin{proof}
    The result follows from a standard approximation argument combined with a priori energy estimates; see \cite{MaJiangWan2016} for the reference.
\end{proof}
By the uniqueness of $H^1$ strong solutions in Theorem~\ref{T2},
the local $H^s$ solution constructed in
Theorem~\ref{T3} coincides with the global $H^1$ solution
on its interval of existence.
To extend this local solution globally in time,
we derive a higher-order a priori estimate using the $H^1$ bounds satisfied by this global solution.
We work on the $\dot{H}^s$ estimate for $(u, v, \te)$, and one has the following Lemma.

\begin{lemma}\label{L5}
Let $s>1$. 
Let $(u, v, \te)$ be the strong solution of \eqref{E11} constructed in Theorem~\ref{T3}.
For each time $t\in[0,T)$, we define $A_s(t)$ and $D_s(t)$ by
\[
A_s (t) := 
\|(\Lam^s u, \Lam^s v, \Lam^s \te) (t)\|_{L^2},
\qquad \text{and}\qquad
D_s (t) :=  \|(\Lam^{s}\nabla u, \Lam^s \nabla v) (t)\|_{L^2}.
\]
Then, there exists a constant $C_s>0$ such that
\begin{equation*}
    \frac{\ud}{\dt}A_s(t)^2 + D_s^2(t) 
    \leq C_s(1+\|v(t)\|_{L^2}^2)\, \|(\nabla u,\nabla v)(t)\|_{L^2}^2 \, A_s^2 
    + \left|
        \int_{\bR^2} \big(\Lambda^s (u \cdot \nabla \te) - u\cdot\Lam^s\nabla\te\big) \Lambda^s \te\,\right|
\end{equation*}
\end{lemma}

\begin{proof}
We multiply $\Lam^{2s} u$, $\Lam^{2s} v$, and $\Lam^{2s} \theta$ on the both sides of $\eqref{E11}_1$-$\eqref{E11}_3$ respectively and integrate over $\mathbb{R}^2$ to deduce that
\begin{equation}\label{E47}
\begin{split}
&\frac{1}{2} \frac{\ud}{\dt}
\|(\Lambda^s u,\Lambda^s v,\Lambda^s \te)\|_{L^2}^2
+ \|(\Lambda^s \nabla u, \Lambda^s \nabla v)\|_{L^2}^2  \\
&= - \int_{\bR^2} \Lambda^s ((u \cdot \nabla) u) \cdot \Lambda^s u
+ \int_{\bR^2} \Lambda^s (v \otimes  v) : \Lambda^s \nabla u 
- \int_{\bR^2} \Lambda^s ((u \cdot \nabla) v) \cdot \Lambda^s v
- \int_{\bR^2} \Lambda^s ((v \cdot \nabla) u) \cdot \Lambda^s v\\
&\quad - \int_{\bR^2} \Lambda^s (u \cdot \nabla \te) \Lambda^s \te\\
& =: \sum_{k=1}^5 \cI_k.
\end{split}
\end{equation}
We use a standard commutator estimate with the divergence-free condition of $u$ to obtain
\begin{align*}
|\cI_1|
&\lesssim \norm{\nabla u}_{L^2}
\norm{\Lam^s u}_{L^4}^2
\lesssim \norm{\nabla u}_{L^2}
\norm{\Lam^s u}_{L^2}
\norm{\Lam^{s + 1} u}_{L^2}
\lesssim \norm{\nabla u}_{L^2} A_s D_s, \\
|\cI_2|
&\lesssim \norm{v}_{L^4}
\norm{\Lam^s v}_{L^4}
\norm{\Lam^{s + 1} u}_{L^2} 
\lesssim \norm{v}_{L^2}^{\frac{1}{2}}
\norm{\nabla v}_{L^2}^{\frac{1}{2}}
\norm{\Lam^s v}_{L^2}^{\frac{1}{2}}
\norm{\Lam^{s + 1} v}_{L^2}^{\frac{1}{2}}
\norm{\Lam^{s + 1} u}_{L^2}
\lesssim 
\norm{v}_{L^2}^{\frac{1}{2}}
\norm{\nabla v}_{L^2}^{\frac12} 
A_s^{\frac12} 
D_s^{\frac32}, \\
|\cI_3|
&\lesssim \norm{\Lam^s u}_{L^4}
\norm{\nabla v}_{L^2}
\norm{\Lam^s v}_{L^4}
+ \norm{\nabla u}_{L^2}
\norm{\Lam^s v}_{L^4}^2 \\
&\lesssim \norm{ \Lam^s u}_{L^2}^{\frac{1}{2}}
\norm{ \Lam^{s + 1} u}_{L^2}^{\frac{1}{2}}
\norm{ \nabla v}_{L^2}
\norm{ \Lam^s v}_{L^2}^{\frac{1}{2}}
\norm{ \Lam^{s + 1} v}_{L^2}^{\frac{1}{2}}
+ \norm{ \Lam^s v}_{L^2}
\norm{ \Lam^{s + 1} v}_{L^2}
\norm{ \nabla u}_{L^2}
\lesssim 
\norm{(\nabla u, \nabla v)}_{L^2} A_s D_s, \\
|\cI_4|
&\lesssim \norm{ \Lam^s v}_{L^4}^2
\norm{ \nabla u}_{L^2}
+ \norm{ \Lam^s v}_{L^4}
\norm{ v}_{L^4}
\norm{\Lam^{s + 1} u}_{L^2} \\
&\lesssim 
\norm{ \Lam^s v}_{L^2}
\norm{ \Lam^{s + 1} v}_{L^2}
\norm{ \nabla u}_{L^2}
+
\norm{v}_{L^2}^{\frac{1}{2}}
\norm{\nabla v}_{L^2}^{\frac{1}{2}}
\norm{\Lam^s v}_{L^2}^{\frac{1}{2}}
\norm{\Lam^{s + 1} v}_{L^2}^{\frac{1}{2}}
\norm{\Lam^{s + 1} u}_{L^2}\\
&
\lesssim 
\norm{\nabla u}_{L^2} A_s D_s
+
\norm{v}_{L^2}^{\frac{1}{2}}
\norm{\nabla v}_{L^2}^{\frac12} 
A_s^{\frac12} 
D_s^{\frac32}.
\end{align*}
Hence, by using Young's inequality, we have
\begin{align*}
\sum_{k=1}^{4} |\cI_k| &
\leq 
C\left( \norm{(\nabla u, \nabla v)}_{L^2} A_s D_s +\norm{v}_{L^2}^{\frac{1}{2}}
\norm{\nabla v}_{L^2}^{\frac12} 
A_s^{\frac12} 
D_s^{\frac32}
\right)\\
& \leq C (1+\|v\|_{L^2}^2)\, 
\norm{(\nabla u, \nabla v)}_{L^2}^2\, A_s^2
+
\frac{1}{2} D_s^2.
\end{align*}

Thanks to divergence-free condition of $u$, the integral $\cI_5$ is understood as
    \begin{equation*}
      \left| \cI_5\right| = \left| \int_{\bR^2} \Lambda^s (u \cdot \nabla \te) \Lambda^s \te \right| 
      =\left|
        \int_{\bR^2} \big(\Lambda^s (u \cdot \nabla \te) - u\cdot\nabla\Lam^s\te\big) \Lambda^s \te \right|. 
    \end{equation*}
   Note that the right-hand side is well defined by the Kato--Ponce
commutator estimate (see \cite{Li2019KP}).
Combining these estimates 
completes the proof of Lemma~\ref{L5}.
\end{proof}

\begin{remark}\label{Remark3}
Lemma~\ref{L5} also holds for $0<s<1$ when $(u, v, \te)$
is a global $H^1$ solution, using the same computation. We shall use this observation later.
\end{remark}

The following corollary extends the $H^1$ result of
Theorem~\ref{T2} to $H^s$ for every $s>1$.
Indeed, Dong--Wu--Ye \cite{DongWuYe2019} established global $H^s$
regularity for \eqref{E11} when $s > 2$, but their statement
does not explicitly cover $1<s\leq2$.
For completeness, we give the details of the continuation
argument for the full range $s>1$.

\begin{corollary}\label{Cor1}
Let $s > 1$ and let the initial data $(u_0,v_0,\theta_0)$ belongs to $H^s(\bR^2)$ with $\dv u_0=0$. Then, the Cauchy problem \eqref{E11} has a unique global solution $(u, v, \te) \in \mathcal{C} ([0,T]; H^s (\bR^2))$ for all $T>0$. 
\end{corollary}

\begin{proof}
Recalling \eqref{E47}, for the term $\cI_5$, we see that
\begin{equation}\label{E36}
\begin{aligned}
|\cI_5|
&= \left| \int_{\bR^2} \Lambda^s (u \cdot \nabla \te) \Lambda^s \te \right|
= \left| \int_{\bR^2} \left(\Lambda^s (u \cdot \nabla \te)
- u\cdot \nabla \Lam^s \te\right) \Lambda^s \te \right| \\
&\lesssim \|\nabla u\|_{L^{\infty}}
\|\Lam^s \te\|_{L^2}^2
+ \begin{cases}
\norm{ \Lam^{s} u}_{L^{\frac{2}{s - 1}}}
\norm{ \nabla \te}_{L^{\frac{2}{2-s}}}
\norm{ \Lam^s \te}_{L^2} 
&\, \text{ if} \ 1 < s < 2, \\
\norm{ \Lam^{s} u}_{L^4}
\norm{ \nabla \te}_{L^4}
\norm{ \Lam^s \te}_{L^2} 
&\,\text{ if} \ s \ge 2,
\end{cases}
\\
&\lesssim \|\nabla u\|_{L^{\infty}}
\|\Lam^s \te\|_{L^2}^2
+
\begin{cases}
\norm{ \Delta u}_{L^2}
\norm{ \Lam^s \te}_{L^2}^2 
&\, \text{ if} \ 1 < s < 2, \\
\left( \norm{ \Delta u}_{L^2}^{\frac{1}{2s - 2}}
\norm{ \Lam^{s + 1} u}_{L^2}^{\frac{2s - 3}{2s - 2}}
\right)
\left( \norm{ \nabla \te}_{L^2}^{\frac{2s - 3}{2s - 2}}
\norm{ \Lam^s \te}_{L^2}^{\frac{1}{2s - 2}}
\right) 
\norm{ \Lam^s \te}_{L^2} 
&\,\text{ if} \ s \ge 2,
\end{cases}
\\
&\lesssim \|\nabla u\|_{L^{\infty}}
\|\Lam^s \te\|_{L^2}^2
+
\begin{cases}
\norm{ \Delta u}_{L^2}
\norm{ \Lam^s \te}_{L^2}^2 
&\, \text{ if} \ 1 < s < 2, \\
\delta\norm{ \Lam^{s + 1} u}_{L^2}^2
+
\norm{ \Delta u}_{L^2}^{\frac{2}{2s - 1}}
\norm{ \nabla \te}_{L^2}^{\frac{4s - 6}{2s - 1}} 
\norm{ \Lam^s \te}_{L^2}^2 
&\,\text{ if} \ s \ge 2,
\end{cases}\\
&\lesssim 
\big(\|\nabla u\|_{L^{\infty}}
+
\norm{ \Delta u}_{L^2}
+
\norm{ \nabla \te}_{L^2}^2
\big)
\|\Lam^s \te\|_{L^2}^2
+
\delta\norm{ \Lam^{s + 1} u}_{L^2}^2.
\end{aligned}
\end{equation}

For $s \ge 2$, the last step uses Young's inequality with $2 / (2s - 1) + (4s - 6) / (2 (2s - 1)) = 1$.

Let $\cE_0:=\|(u_0,v_0,\te_0)\|_{L^2}^2$. 
After taking $\delta>0$ sufficiently small and using the absorbing scheme, we combine this with Lemma~\ref{L5}. Together with Lemma~\ref{Lem:L2}, applying Gr\"onwall's inequality to the resulting estimate yields
\begin{equation}\label{eq:monotone}
    \begin{aligned}
        A_s(t)^2 +\int_0^t D_s(\tau)^2 \dtau 
        &\leq  2A_s(0)^2 \exp\left(C\int_0^t (1+\|v(\tau)\|_{L^2}^2)
\|(\nabla u,\nabla v)(\tau)\|_{L^2}^2
        +\|\nabla u(\tau)\|_{L^{\infty}}
+\norm{ \Delta u(\tau)}_{L^2}
+
\norm{ \nabla \te(\tau)}_{L^2}^2
\dtau\right)\\
        &\leq 2 A_s(0)^2 \exp\left(C (\cE_0+\cE_0^2)
+ C\int_0^t(\|\nabla u(\tau)\|_{L^{\infty}}
+\norm{ \Delta u(\tau)}_{L^2}
+
\norm{ \nabla \te(\tau)}_{L^2}^2)\,\dtau\right).
    \end{aligned}
    \end{equation}

By Theorem~\ref{T2},  for every $T>0$, we have
\[
\begin{aligned}
    &\int_0^T
    \left(
        \|\nabla u(\tau)\|_{L^\infty}
        +\|\Delta u(\tau)\|_{L^2}
        +\|\nabla\te(\tau)\|_{L^2}^2
    \right)\,\ud\tau
\\
    &\leq
    \|\nabla u\|_{L^1([0,T];L^\infty)}
    +T^{1/2}\|\Delta u\|_{L^2([0,T];L^2)}
    +T\|\nabla\te\|_{L^\infty([0,T];L^2)}^2
    <\infty.
\end{aligned}
\]
This excludes finite-time blow-up of the $H^s$ norm.
The local-in-time solution constructed in Theorem~\ref{T3} therefore extends globally by a standard continuation argument. The proof is finished.
\end{proof}

\subsection{Proof of main theorem}
As in the $H^1$ stability analysis, we consider the following \textit{cross term} involving $v$ and $\te$.
We multiply $\Lam^{2s - 2} \nabla \theta$ on $\eqref{E11}_2$ and $-\Lam^{2s - 2} \dv v$ on $\eqref{E11}_3$. 
Then, we integrate over $\mathbb{R}^2$ to obtain 
\begin{align*}
\int_{\bR^2} v_t \cdot \Lam^{2s - 2} \nabla \theta
+ \int_{\bR^2} (u \cdot \nabla) v \cdot \Lam^{2s - 2} \nabla \theta
+ \int_{\bR^2} (v \cdot \nabla) u \cdot \Lam^{2s - 2} \nabla \theta
- \int_{\bR^2} \Delta v \cdot \Lam^{2s - 2} \nabla \theta
&= \int_{\bR^2} | \Lam^{s - 1} \nabla \theta|^2 \\
\int_{\bR^2} v \cdot \Lam^{2s - 2} \nabla \theta_t 
- \int_{\bR^2} (u \cdot \nabla \theta)\Lam^{2s - 2} \dv v 
&= - \int_{\bR^2} |\Lam^{s - 1} \dv v|^2.
\end{align*}
Adding these identities and using the integration by parts, we obtain
\begin{align*}
\int_{\bR^2} |\Lam^s \theta|^2
&+ \int_{\bR^2} \Lam^{s - 1} \Delta v \cdot \Lam^{s - 1} \nabla \theta
- \int_{\bR^2} |\Lam^{s - 1} \dv v|^2 
- \frac{\ud}{\dt} \int_{\bR^2} \Lam^{s - 1} v \cdot \nabla \Lam^{s - 1} \theta \\
&= \int_{\bR^2} \Lam^{s - 1} ((u \cdot \nabla) v) \cdot \Lam^{s - 1} \nabla \theta
+ \int_{\bR^2} \Lam^{s - 1} ((v \cdot \nabla) u) \cdot \Lam^{s - 1} \nabla \theta
- \int_{\bR^2} \Lam^{s - 1} (u \cdot \nabla \theta) \Lam^{s - 1} \dv v.
\end{align*}
Young's inequality reveals that
\begin{equation}\label{E48}
\begin{split}
&\frac{1}{2} \|\Lam^s \te\|_{L^2}^2
- \|\Lam^s v\|_{H^1}^2
- \frac{\ud}{\dt} \int_{\bR^2} \Lam^{s - 1} v \cdot \Lam^{s - 1} \nabla \te \\
&\quad \le \int_{\bR^2} \Lam^{s - 1} ((u \cdot \nabla) v) \cdot \Lam^{s - 1} \nabla \theta
+ \int_{\bR^2} \Lam^{s - 1} ((v \cdot \nabla) u) \cdot \Lam^{s - 1} \nabla \theta
- \int_{\bR^2} \Lam^{s - 1} (u \cdot \nabla \theta) \Lam^{s - 1} \dv v.
\end{split}
\end{equation}
Multiplying the inequality~\eqref{E48} by $1/8$ and adding the energy identities~\eqref{E47} at orders $s$ and $s-1$, we arrive at
\begin{equation}\label{E61}
\begin{aligned}
&\frac{1}{2} \frac{\ud}{\dt}
\left( \|\Lambda^s u\|_{L^2}^2
+ \|\Lambda^s v\|_{L^2}^2
+ \|\Lambda^s \te\|_{L^2}^2 
+\|\Lambda^{s-1} u\|_{L^2}^2
+ \|\Lambda^{s-1} v\|_{L^2}^2
+ \|\Lambda^{s-1} \te\|_{L^2}^2 
- \frac{1}{4} \int_{\bR^2} \Lam^{s - 1} v \cdot \Lam^{s - 1} \nabla \theta
\right) \\
&\quad + 
\left(\|\Lambda^{s+1} u\|_{L^2}^2 
+ \frac{1}{2} \|\Lambda^{s+1}  v\|_{L^2}^2
+\frac{1}{16} \| \Lam^{s} \theta \|_{L^2}^2
+ \| \Lam^{s} u\|_{L^2}^2
+ \frac{1}{2} \| \Lam^{s} v \|_{L^2}^2
\right) \\
&\le - \int_{\bR^2} \Lambda^s ((u \cdot \nabla) u) \cdot \Lambda^s u
+ \int_{\bR^2} \Lambda^s (v \otimes  v) : \Lambda^s \nabla u 
- \int_{\bR^2} \Lambda^s ((u \cdot \nabla) v) \cdot \Lambda^s v
- \int_{\bR^2} \Lambda^s ((v \cdot \nabla) u) \cdot \Lambda^s v\\
&\quad - \int_{\bR^2} \Lambda^s (u \cdot \nabla \te)
\Lambda^s \te
+ \frac{1}{8} \int_{\bR^2} \Lam^{s - 1} (u \cdot \nabla) v \cdot \Lam^{s - 1} \nabla \theta \\
&\quad + \frac{1}{8} \int_{\bR^2} \Lam^{s - 1} (v \cdot \nabla) u \cdot \Lam^{s - 1} \nabla \theta
- \frac{1}{8} \int_{\bR^2} \Lam^{s - 1} (u \cdot \nabla \theta) \Lam^{s - 1} \dv v\\
&\quad - \int_{\bR^2} \Lambda^{s-1} ((u \cdot \nabla) u) \cdot \Lambda^{s-1} u
+ \int_{\bR^2} \Lambda^{s-1} (v \otimes  v) : \Lambda^{s-1} \nabla u 
- \int_{\bR^2} \Lambda^{s-1} ((u \cdot \nabla) v) \cdot \Lambda^{s-1} v \\
&\quad - \int_{\bR^2} \Lambda^{s-1} ((v \cdot \nabla) u) \cdot \Lambda^{s-1} v
- \int_{\bR^2} \Lambda^{s - 1} (u \cdot \nabla \te) \Lambda^{s - 1} \te \\
&\quad =: \sum_{k=1}^{13} \cI_k.
\end{aligned}
\end{equation}
Hence, it is natural to define
\begin{equation}\label{def:X_s}
    X_s(t):= \left( \|\Lambda^s u\|_{L^2}^2
+ \|\Lambda^s v\|_{L^2}^2
+ \|\Lambda^s \te\|_{L^2}^2 
+\|\Lambda^{s-1} u\|_{L^2}^2
+ \|\Lambda^{s-1} v\|_{L^2}^2
+ \|\Lambda^{s-1} \te\|_{L^2}^2 
- \frac{1}{4} \int_{\bR^2} \Lam^{s - 1} v \cdot \Lam^{s - 1} \nabla \theta
\right)^{\frac{1}{2}}(t),
\end{equation}
and
\begin{equation}\label{def:Y_s}
  Y_s(t):= \left(\|\Lambda^{s+1} u\|_{L^2}^2 
+ \frac{1}{2} \|\Lambda^{s+1}  v\|_{L^2}^2
+\frac{1}{16} \| \Lam^{s} \theta \|_{L^2}^2
+ \| \Lam^{s} u\|_{L^2}^2
+ \frac{1}{2} \| \Lam^{s} v \|_{L^2}^2
\right)^{\frac{1}{2}}(t).
\end{equation}

\begin{lemma}\label{Lem:Hs_key}
Let $s > 1$ and $1 < p < 2$. Let $(u, v, \te)$ be the $H^s$ solution constructed in Corollary~\ref{Cor1}, and let $\cE_0 := \|(u_0,v_0,\te_0)\|_{L^2}^2$. For $X_s$ and $Y_s$ defined in \eqref{def:X_s}--\eqref{def:Y_s}, there exists $C = C_{s,p} > 0$ such that
\begin{equation*}
X_s (t)^2 + \int_r^t Y_s(\tau)^2 d \tau
\le 2X_s(r)^2
\exp\left(
C(\cE_0+\cE_0^2)+
C\int_r^t 
\norm{\nabla u (\tau)}_{L^\infty}^p
+
\norm{\Del u (\tau)}_{L^2}^2
+
\norm{\nabla \theta (\tau)}_{L^2}^2\, \dtau
\right),
\end{equation*}
for every $0\leq r<t$.
\end{lemma}

\begin{proof}
We start from \eqref{E61}.
The terms $\cI_1,\ldots,\cI_4$ and
$\cI_9,\ldots,\cI_{12}$ can be estimated as in the proof
of Lemma~\ref{L5} (see also Remark~\ref{Remark3}), yielding
    \begin{equation*}
    \begin{aligned}
         \sum_{k=1}^4|\cI_k|
    +\sum_{k=9}^{12}|\cI_k|
    &\leq
C_s\|\nabla(u,v)\|_{L^2}X_sY_s
+ C_s\|v\|_{L^2}^{1/2}\|\nabla v\|_{L^2}^{1/2}X_s^{1/2}Y_s^{3/2} \\
    &\leq 
C_{s,\delta}\bigl(1+\|v\|_{L^2}^2\bigr)
    \|\nabla(u,v)\|_{L^2}^2 X_s^2
    +\delta Y_s^2.
    \end{aligned}
    \end{equation*}
To estimate $\cI_5$, from~\eqref{E36}, we have
\begin{align*}
|\cI_5|
&\leq C\big(\|\nabla u\|_{L^{\infty}}
+
\norm{ \Delta u}_{L^2}
+
\norm{ \nabla \te}_{L^2}^2
\big)
\|\Lam^s \te\|_{L^2}^2
+
\delta\norm{ \Lam^{s + 1} u}_{L^2}^2 \\
& \leq C\big(\|\nabla u\|_{L^{\infty}}^p
+
\norm{ \Delta u}_{L^2}^2
+
\norm{ \nabla \te}_{L^2}^2
\big)
\|\Lam^s \te\|_{L^2}^2
+
\delta\norm{ \Lam^{s} \te}_{L^2}^2
+
\delta\norm{ \Lam^{s + 1} u}_{L^2}^2
\\
&\leq
C_{s,p,\delta} \big(\|\nabla u\|_{L^{\infty}}^p
+
\norm{ \Delta u}_{L^2}^2
+
\norm{ \nabla \te}_{L^2}^2
\big)
X_s^2
+
\delta Y_s^2.
\end{align*}
In the second line, we used the inequalities $a \leq \delta +C(p,\delta)a^p$ and $b \leq \delta +\frac{1}{4\delta}b^2$ for $a,b\geq0$.

For the term $\cI_6$, we use 
$\|w\|_{L^{\infty}} \lesssim \|w\|_{L^2}^{1/2}\|\Delta w\|_{L^2}^{1/2}$ (see \cite[Lemma A.1]{MR5029968}) to obtain
\begin{align*}
|\cI_6|
&\lesssim \left| \int_{\bR^2} \Lambda^{s - 1} (u \cdot \nabla v) \cdot \Lambda^{s - 1} \nabla \te \right| \\
&\lesssim \|u\|_{L^{4}}
\|\Lam^s v\|_{L^4}
\|\Lam^s \te\|_{L^2}
+
\norm{ \Lam^{s-1} u}_{L^{\infty}}
\norm{ \nabla v}_{L^2}
\norm{ \Lam^s \te}_{L^2} \\
&\lesssim \|u\|_{L^2}^{\frac{1}{2}}
\|\nabla u\|_{L^2}^{\frac{1}{2}}
\|\Lam^s v\|_{L^2}^{\frac{1}{2}}
\|\Lam^{s+1} v\|_{L^2}^{\frac{1}{2}}
\|\Lam^s \te\|_{L^2}
+
\norm{ \Lam^{s-1} u}_{L^2}^{\frac{1}{2}}
\norm{ \Lam^{s+1} u}_{L^2}^{\frac{1}{2}}
\norm{ \nabla v}_{L^2}
\norm{ \Lam^s \te}_{L^2}^\frac{1}{2}
\norm{ \Lam^s \te}_{L^2}^\frac{1}{2}\\
& \lesssim \|u\|_{L^2}^{\frac{1}{2}}
\|\nabla u\|_{L^2}^{\frac{1}{2}}
X_s^{\frac{1}{2}}
Y_s^{\frac{3}{2}}
+
\norm{ \nabla v}_{L^2}
X_s
Y_s \\
&\lesssim
\left(
\|u\|_{L^2}^2
\|\nabla u\|_{L^2}^2
+
\norm{ \nabla v}_{L^2}^2
\right)
X_s^2
+
\delta
Y_s^2
\end{align*}
Similar computation yields that
\begin{align*}
&|\cI_7|
\lesssim \left| \int_{\bR^2} \Lambda^{s - 1} (v \cdot \nabla u) \cdot \Lambda^{s - 1} \nabla \te \right| 
\lesssim 
\|v\|_{L^2}^{\frac12}
\|\nabla v\|_{L^2}^{\frac12}
X_s^{\frac12}
Y_s^{\frac32}
+
\norm{ \nabla u}_{L^2}
X_s
Y_s
\lesssim 
\left(
\|v\|_{L^2}^2
\|\nabla v\|_{L^2}^2
+
\norm{ \nabla u}_{L^2}^2
\right)
X_s^2
+
\delta
Y_s^2,
\\
&|\cI_8|
\lesssim 
\left| \int_{\bR^2} \Lam^{s - 1} (u \cdot \nabla \theta) \Lam^{s - 1} \dv v \right| 
\lesssim 
\|u\|_{L^2}^{\frac12}
\|\nabla u\|_{L^2}^{\frac12}
X_s^{\frac12}
Y_s^{\frac32}
+
\norm{ \nabla \te}_{L^2}
X_s
Y_s
\lesssim 
\left(
\|u\|_{L^2}^2
\|\nabla u\|_{L^2}^2
+
\norm{ \nabla \te}_{L^2}^2
\right)
X_s^2
+
\delta
Y_s^2.
\end{align*}
Lastly,
\begin{align*}
    |\cI_{13}| &= \left|\int_{\bR^2} 
    \left(\Lambda^{s-1} (u \cdot \nabla \te)- u 
    \cdot 
    \nabla \Lam^{s-1}\te\right) \Lambda^{s-1} \te
    \right|\\
&\lesssim \|\nabla u\|_{L^2}
\|\Lam^{s-1} \te\|_{L^4}^2
+ 
\norm{ \Lam^{s-1} u}_{L^4}
\norm{ \nabla \te}_{L^2}
\norm{ \Lam^{s-1} \te}_{L^4}\\
&\lesssim \|\nabla u\|_{L^2}
\|\Lam^{s-1} \te\|_{L^2}
\|\Lam^{s} \te\|_{L^2}
+ 
\norm{ \Lam^{s-1}u}_{L^2}^{\frac12}
\norm{ \Lam^{s}u}_{L^2}^{\frac12}
\norm{ \nabla \te}_{L^2}
\norm{ \Lam^{s-1}\te}_{L^2}^{\frac12}
\norm{ \Lam^{s} \te}_{L^2}^{\frac12}\\
&\lesssim 
\left(\norm{\nabla u}_{L^2}
+\|\nabla \te\|_{L^2}\right)
X_s Y_s\\
&\lesssim 
\left(\norm{\nabla u}_{L^2}^2
+\|\nabla \te\|_{L^2}^2\right)
X_s^2
+
\delta
Y_s^2.
\end{align*}

Combining these, we deduce that
\begin{equation*}
\begin{split}
    \frac{d}{dt}X_s(t)^2 + Y_s(t)^2 
    \leq 
    C_{s,p,\delta}
    \left(
    \bigl(1+\|(u,v)\|_{L^2}\bigr)^2
    \|\nabla(u,v)\|_{L^2}^2
    +\|\nabla u\|_{L^{\infty}}^p
    +\norm{ \Delta u}_{L^2}^2
    +\|\nabla \te\|_{L^2}^2
    \right)\cdot X_s^2
    +
    C_s\,\delta Y_s^2 
\end{split}
\end{equation*}
By choosing $\delta>0$ such that $C_s\delta = \frac{1}{2}$, we absorb $C_s\delta Y_s^2$ into left side. Applying Gr\"onwall inequality completes the proof.
\end{proof}

\begin{proof}[Proof of Theorem~\ref{Thm1}]
We fix $p=3/2$ and choose the constant $0< \ep_0 \leq \ep_1$, where $\ep_1$ is the constant in
Proposition~\ref{prop_stab_H1}.
The case $s=1$ is Proposition~\ref{prop_stab_H1}, hence we assume $s>1$.

\noindent
\emph{(i) \textbf{Stability.}}

First, we establish \eqref{E10}. 
By taking $r=0$ in Lemma~\ref{Lem:Hs_key}, for all $t \geq 0$, we have
\begin{equation}\label{E69}
\begin{aligned}
X_s (t)^2 + \int_0^t Y_s(\tau)^2 d \tau
&
\le 2X_s(0)^2
\exp\left(
C(\cE_0+\cE_0^2)+
C\int_0^{\infty} 
\norm{\nabla u (\tau)}_{L^\infty}^p
+
\norm{\Del u (\tau)}_{L^2}^2
+
\norm{\nabla \theta (\tau)}_{L^2}^2\, \dtau
\right)\\
&
\le 2X_s(0)^2
\exp\left(
C(\cE_0+\cE_0^2)+
C\|(u_0,v_0,\te_0)\|_{H^1}^p
+ C\|(u_0,v_0,\te_0)\|_{H^1}^2
\right)\\
&\le
C X_s(0)^2,
\end{aligned}
\end{equation}
where the constant $C$ depends only on $s>1$. 
In \eqref{E69}, the second inequality follows from \eqref{E12} and \eqref{E53}, where the latter is obtained from Lemma~\ref{Prop1} under the smallness assumption on the initial $H^1$ norm.
The last inequality uses $\|(u_0,v_0,\te_0)\|_{H^1}<\ep_0\leq1$. 
Combining \eqref{E69} with Lemma~\ref{Lem:L2} and Proposition~\ref{prop_stab_H1} proves \eqref{E10}.

\vspace{1em}
\noindent
\emph{(ii) \textbf{Temporal decay.}}

We divide the proof into three steps.
We recall that
\[
A_s (t) := 
\|(\Lam^s u, \Lam^s v, \Lam^s \te) (t)\|_{L^2},
\qquad \text{and}\qquad
D_s (t) :=  \|(\Lam^{s}\nabla u, \Lam^s \nabla v) (t)\|_{L^2}.
\]

\medskip
\noindent
\emph{Step 1) Monotonicity of the $\dot H^s$ norm on a dyadic interval.}

First, we claim that the following: For every $T \geq1$ and $T \leq r \leq t \leq 2T$, assuming the smallness of the $H^1$ norm of initial data, we have
\begin{equation}\label{claim:mono}
    \|\Lam^s (u, v, \te)(t)\|_{L^2}\leq    C \|\Lam^s (u, v, \te)(r)\|_{L^2},
\end{equation}
where the constant $C$ depends only on $s>1$. Indeed, by just changing $0$ to $r$ in \eqref{eq:monotone}, we have
\begin{equation*}
    \begin{aligned}
        A_s(t)^2 +\int_r^{t} D_s(\tau)^2 \dtau
        &\leq 2 A_s(r)^2 \exp\left(C (\cE_0+\cE_0^2)
+ C\int_r^{t}(\|\nabla u(\tau)\|_{L^{\infty}}
+\norm{ \Delta u(\tau)}_{L^2}
+
\norm{ \nabla \te(\tau)}_{L^2}^2)\,\dtau\right)\\
&\leq 2 A_s(r)^2 \exp\left(C (\cE_0+\cE_0^2)
+ C\int_T^{2T}(\|\nabla u(\tau)\|_{L^{\infty}}
+\norm{ \Delta u(\tau)}_{L^2}
+
\norm{ \nabla \te(\tau)}_{L^2}^2)\,\dtau\right).
    \end{aligned}
    \end{equation*}
This estimate does not require any smallness assumption.
Since $[r,t]\subset[T,2T]$, Lemma~\ref{lem:dyadic_int} and \eqref{E12} yield, under the assumed smallness condition,
\begin{equation}\label{eq:useful1605}
    \begin{aligned}
        A_s(t)^2 +\int_r^{t} D_s(\tau)^2 \dtau
        &\leq 2 A_s(r)^2 \exp\left(
C_s\|(u_0,v_0,\te_0)\|_{H^1}\right) 
 \leq C_s\, A_s(r)^2.
    \end{aligned}
    \end{equation}
This proves the claim.

\medskip
\noindent
\emph{Step 2) Minimum selection argument.}

Fix $T\ge1$. Since the map $t\mapsto\|\Lam^s(u, v, \te)(t)\|_{L^2}^2$ is continuous, choose a minimizer $r^*\in[T,3T/2]$.
By Lemma~\ref{Lem:Hs_key} with $r = T$ and $t = 3T/2$, as in \eqref{E69}, we can show that
\[
X_s (3T/2)^2
+ \int_T^{3T/2} Y_s (\tau)^2 \dtau
\le C X_s (T)^2,
\]
where the constant $C$ depends only on $s$.
Thus,
\begin{align*}
\norm{\Lam^s  (u, v, \te) (r^*)}_{L^2}^2
&\le \frac{2}{T} \int_{T}^{3T/2}
\norm{\Lam^s  (u, v, \te) (\tau)}_{L^2}^2 \dtau \\
&\le \frac{C}{T} \int_{T}^{3T/2}
Y_s(\tau)^2 \dtau
\leq \frac{C}{T} X_s(T)^2
\leq \frac{C}{T} 
\left(\norm{\Lam^s  (u, v, \te) (T)}_{L^2}^2
+ \norm{\Lam^{s-1}  (u, v, \te) (T)}_{L^2}^2
\right). 
\end{align*}
Using \eqref{claim:mono} with $t=2T$ and $r=r^*$, we have
\begin{equation}\label{eq:iter1449}
\norm{\Lam^s  (u, v, \te) (2T)}_{L^2}^2 
\leq C
\norm{\Lam^s  (u, v, \te) (r^*)}_{L^2}^2
\leq
\frac{C}{T} 
\left(\norm{\Lam^s  (u, v, \te) (T)}_{L^2}^2
+
\norm{\Lam^{s-1}  (u, v, \te) (T)}_{L^2}^2
\right). 
\end{equation}

\medskip
\noindent
\emph{Step 3) Induction argument and iteration.}

If $0 < \gamma < 1$, then it follows from Proposition~\ref{prop_stab_H1}, Lemma \ref{Lem:L2}, and the interpolation inequality that
\[
\|(\Lambda^\gamma u, \Lambda^\gamma v, \Lambda^\gamma \te)(t)\|_{L^2}
\leq
\|(\nabla u, \nabla v, \nabla \te)(t)\|_{L^2}^{\gamma}
\|(u, v, \te)(t)\|_{L^2}^{1 - \gamma}
\leq \left(\cM_1 (t) (1 + t)^{- 1/2}
\right)^\gamma
\|(u_0,v_0,\te_0)\|_{L^2}^{1 - \gamma}
\]
and
\begin{align*}
\left( \int_t^\infty \norm{(\Lambda^\gamma \nabla u, \Lambda^\gamma \nabla v) (t)}_{L^2}^2 dt \right)^{1/2}
&\le \left( \int_t^\infty \norm{(\Delta u, \Delta v) (t)}_{L^2}^2 dt \right)^{\gamma/2}
\left( \int_t^\infty \norm{(\nabla u, \nabla v) (t)}_{L^2}^2 dt \right)^{(1 - \gamma) /2} \\
&\le \left(\cM_1 (t) (1 + t)^{- 1/2}
\right)^\gamma
\|(u_0,v_0,\te_0)\|_{L^2}^{1 - \gamma}.
\end{align*}
Thus \eqref{E13} holds at order $\gamma$ with
\[
\cM_{\gamma} (t)
:= \cM_1 (t)^\gamma
\|(u_0,v_0,\te_0)\|_{L^2}^{1 - \gamma}.
\]
The case $\gamma = 1$ is Proposition~\ref{prop_stab_H1}.

Hence, for the induction argument, suppose that \eqref{E13} holds with $s$ replaced by $s-1$.
In the following induction, $\cM_{s-1}$ denotes a constant
multiple of the square of the preceding majorant.
Below, we set
\begin{equation*}
\varphi(T):=\norm{\Lam^s  (u, v, \te) (T)}_{L^2}^2. 
\end{equation*}
The induction hypothesis and \eqref{eq:iter1449} yield
\begin{equation*}
\vphi (2T) \le \frac{C}{T} \vphi (T) 
+ \frac{\cM_{s - 1} (T)}{(1 + T)^s}.
\end{equation*}
where $\cM_{s-1}$ is non-increasing and satisfies
\begin{equation*}
    \sup_{t\geq0} \cM_{s-1}(t) \leq C_{s-1} \|(u_0,v_0,\te_0)\|_{H^{\max\{1,s-1\}}}^2,\quad \text{ and } \quad
    \lim_{t \to \infty}\cM_{s-1}(t) = 0.
\end{equation*}

Iterating the preceding inequality $N$-times, we obtain
\[
\begin{aligned}
\vphi(2^NT)
&\le
\frac{C(N,s)}{(1+T)^N}\vphi(T)
+\frac{C(N,s)\cM_{s-1}(T)}{(1+T)^s}
\leq 
\frac{C(N,s)}{(1+T)^s}\frac{\|(u_0,v_0,\te_0)\|_{H^s}^2}{(1+T)^{N-s}}
+\frac{C(N,s)\cM_{s-1}(T)}{(1+T)^s}.
\end{aligned}
\]
Before adjusting constants, the coefficient of $\vphi (T)$ after $N$ iterations is $C_s^N 2^{- N (N - 1) / 2} T^{- N}$.
We fix $N = [s] + 1$, where $[s]$ denotes the integer part of $s$; thus $N > s$.
Setting $t=2^NT$ and adjusting the constant gives
\[
(1+t)^s\vphi(t)
\le C_s
\left(
    \frac{\|(u_0,v_0,\te_0)\|_{H^s}^2}{(1+t)^{N-s}}
    +\cM_{s-1}(2^{-N}t)
\right),
\qquad t\ge2^N.
\]
We therefore define
\begin{equation*}
    \cM_s(t) := C_s\left(\frac{\|(u_0,v_0,\te_0)\|_{H^s}^2}{(1+t)^{[s]-s+1}} + \cM_{s-1}(2^{-[s]-1}t)\right).
\end{equation*}

Hence, $\cM_s$ is non-increasing and satisfies
\[
\norm{\Lam^s  (u, v, \te) (t)}_{L^2}^2 \le\frac{\cM_s(t)}{(1+t)^s},
\qquad
\sup_{t \ge 0}\cM_s(t)\le C_s\|(u_0,v_0,\te_0)\|_{H^s}^2,
\qquad
\lim_{t\to\infty}\cM_s(t)=0.
\]
Starting from the range $0<s\le1$, this argument applies successively to every $s>1$. 
By \eqref{eq:useful1605}, for $t\ge1$, we also obtain
\begin{equation*}
\int_{t}^{2t}\|\Lam^{s+1}(u,v)(\tau)\|_{L^2}^2\dtau \leq C_s A_s(t)^2
    \leq C_s\frac{\cM_s(t)}{(1+t)^s}.
\end{equation*}
Consequently, by the monotonicity of $\cM_s$,
\[
\begin{aligned}
\int_t^\infty
\|\Lambda^{s+1}(u,v)(\tau)\|_{L^2}^2\,\ud\tau
&=
\sum_{k=0}^\infty
\int_{2^kt}^{2^{k+1}t}
\|\Lambda^{s+1}(u,v)(\tau)\|_{L^2}^2\,\ud\tau\\
&\le
C_s\sum_{k=0}^\infty
\frac{\cM_s(2^kt)}{(1+2^kt)^s}
\le
\frac{C_s
\cM_s(t)}{(1+t)^s},
\qquad \textrm{for }t\ge1.
\end{aligned}
\]
Here $1 + 2^k t \ge 2^{k - 1} (1 + t)$ for $t \ge 1$.
This completes the proof of Theorem~\ref{Thm1}.
\end{proof}

\section{Sharpness of temporal decay}
\label{S4}
It is classical that planar Navier--Stokes flows with radial
vorticity also solve the heat equation, since the nonlinear term
is absorbed by the pressure; see
\cite[p.~224]{Brandolese2004}.
The following lemma gives a related symmetry reduction
for \eqref{E11}.

\begin{lemma}
\label{lem:S4}
Let $s > 2$, and let $(v_0,\te_0)\in H^s(\bR^2)$
be initial data of the form
\[
v_0(x)=a_0(|x|)\frac{x}{|x|},\qquad
\te_0(x)=\Theta_0(|x|),
\]
for some scalar functions $a_0, \Theta_0 : [0,\infty) \to \bR$.
Let $(u, v, \te)\in C([0, \infty) ; H^s(\bR^2))$
be the unique global strong solution of \eqref{E11}
given by Corollary~\ref{Cor1}, with initial data
\begin{equation}\label{eq:sp_initial}
(u, v, \te)\big|_{t=0}=(0,v_0,\te_0).
\end{equation}
Then, $u(t)\equiv0$ for every $t \ge 0$.
Moreover, there exist scalar functions
$a(t, r)$ and $\Theta(t, r)$ such that
\[
v(t, x)=a(t, |x|)\frac{x}{|x|},\qquad
\te(t, x)=\Theta(t, |x|).
\]
Here, the expressions
$a_0(|x|)\frac{x}{|x|}$ and $a(t, |x|)\frac{x}{|x|}$ are understood by continuous
extension at the origin, with $a_0(0)=a(t,0)=0$.
\end{lemma}
\begin{proof}
For any orthogonal matrix $Q$, define
\[
\begin{aligned}
u^Q(t, x)&:=Q^Tu(t,Qx),&
v^Q(t, x)&:=Q^Tv(t,Qx),\\
\te^Q(t, x)&:=\te(t,Qx),&
p^Q(t, x)&:=p(t,Qx).
\end{aligned}
\]
A direct computation shows that
$(u^Q,v^Q,\te^Q,p^Q)$ also solves \eqref{E11}.
Since
\[
v_0(Qx)=Qv_0(x),\qquad
\te_0(Qx)=\te_0(x),
\]
the transformed solution has the same initial data
\eqref{eq:sp_initial}.
Uniqueness therefore gives
\[
(u^Q,v^Q,\te^Q)=(u, v, \te).
\]
Hence, for every orthogonal matrix $Q$,
\[
u(t,Qx)=Qu(t, x),\qquad
v(t,Qx)=Qv(t, x),\qquad
\te(t,Qx)=\te(t, x).
\]

Recall that $e_1=(1,0)$ and $e_2=(0,1)$ are the
standard Cartesian basis vectors, and write
\[
u=u_1e_1+u_2e_2,\qquad
v=v_1e_1+v_2e_2.
\]
For any $x \ne 0$, set $r := |x|>0$. Taking $Q=\operatorname{diag}(1,-1)$ at $x=re_1$ gives
\[
u_2(t,re_1)=v_2(t,re_1)=0.
\]
Define
\[
b(t, r):=u_1(t,re_1),\qquad
a(t, r):=v_1(t,re_1),\qquad
\Theta(t, r):=\te(t,re_1).
\]
We choose a rotation $R$
such that $Re_1=x/r$. Then
\[
\begin{aligned}
u(t, x)&=R u(t,re_1)=b(t, r)e_r,\\
v(t, x)&=R v(t,re_1)=a(t, r)e_r,\\
\te(t, x)&=\te(t,re_1)=\Theta(t, r).
\end{aligned}
\]
This proves the preservation of the stated symmetries.

Since $s > 2$, Sobolev embedding implies that $u(t)$
is continuously differentiable and bounded.
Thus $\dv u=0$ yields
\[
0=\partial_r b+\frac{b}{r}
=\frac1r\partial_r(rb)
\qquad(r>0).
\]
Hence $b(t, r)=c(t)/r$ for $r>0$.
Boundedness near the origin forces $c(t)=0$,
so $u(t)\equiv0$ for every $t \ge 0$.
\end{proof}

\begin{remark}
Let $w=(w_1,w_2)\in C^1(\bR^2)$ satisfy $\dv w=0$ and
\[
w(Qx)=Qw(x)
\qquad\text{for any orthogonal matrix }Q.
\]
The reflection and rotation argument above gives
\[
w(x)=\beta(r)e_r,\qquad r = |x|>0.
\]
Hence
\[
0=\dv w=\frac1r\partial_r(r\beta),
\qquad
\beta(r)=\frac Cr.
\]
Since $w$ is bounded near the origin, necessarily $C=0$.
Continuity then extends $w\equiv0$ on $\bR^2$.

Thus, when considering radially symmetric solutions
of \eqref{E11} in this class, the condition $\dv u_0=0$
already forces $u_0=0$.
This explains the choice of initial barotropic velocity
in Lemma~\ref{lem:S4}.
\end{remark}
\begin{remark}
For the solution in Lemma~\ref{lem:S4}, we have $u\equiv0$,
so the first equation of \eqref{E11} reduces to
\[
\nabla p+\dv(v\otimes v)=0.
\]
We verify this identity directly.
Since $v(t,0)=0$, the regularity of $v$ gives 
$a(t, r)=O(r)$ as $r\downarrow0$.
Below, we set
\[
p(t, x)
:=-a(t, r)^2-\int_0^r\frac{a(t,\rho)^2}{\rho}\,\ud\rho,
\qquad r = |x|.
\]
Then, $p(t,\cdot)\in C^1(\bR^2)$.
For $x \ne 0$, a direct computation gives
\[
\nabla p
=-\left(2a\partial_r a+\frac{a^2}{r}\right)e_r
=-\dv(v\otimes v).
\]
By continuity, this identity extends to $\bR^2$.
\end{remark}

\begin{proof}[Proof of Theorem~\ref{Thm2:sharp}]

For a related annular construction, see \cite[Example~6.3]{Brandolese2016}.
Fix $\ep$ and $h$ as in the theorem. 
Choose $\psi\in\mathcal S(\bR^2)$
with $\sF\psi$ real and radial, such that
\[
\supp(\sF\psi)\subset\{1<|\xi|<2\},\qquad
\|\psi\|_{L^2}=1,
\]
and choose
\[
t_1\ge64,\qquad t_{n+1}>16t_n,\qquad h(t_n)\le4^{-n}.
\]
Set
\begin{equation}\label{E44}
\begin{aligned}
\te_0&:=\sum_{n=1}^\infty\psi_n,
&\psi_n(x)&:=\ep2^{-n}t_n^{-1/2}\psi(x/\sqrt{t_n}),\\
\sF\psi_n(\xi)&=\ep2^{-n}\sqrt{t_n}\,
\sF\psi(\sqrt{t_n}\xi),
&\supp(\sF\psi_n)&\subset
A_n:=\{t_n^{-1/2}<|\xi|<2t_n^{-1/2}\}.
\end{aligned}
\end{equation}
The annuli $A_n$ are pairwise disjoint, and all $A_n$ lie in $\{|\xi|<1/4\}$.
A direct computation shows that
\[
\|\psi_n\|_{H^1}^2
=\int_{\mathbb R^2}(1+|\xi|^2)
|\widehat{\psi_n}(\xi)|^2\,d\xi
\le\left(1+\frac4{t_n}\right)\|\psi_n\|_{L^2}^2
\le\epsilon^24^{-n}\left(\frac{17}{16}\right),
\]
and
\[
\|\theta_0\|_{H^1}^2
=\sum_{n=1}^{\infty}\|\psi_n\|_{H^1}^2
\le\epsilon^2\left(\frac{17}{16}\right)
\sum_{n=1}^{\infty}4^{-n}
\leq \ep^2.
\]
Moreover, we have $\te_0 \in\bigcap_{m\ge1}H^m(\bR^2)$.
Let $(u, v, \te)$ be the unique global strong solution
of \eqref{E11} with initial data $(0,0,\te_0)$,
given by Corollary~\ref{Cor1}.
By Lemma~\ref{lem:S4}, we have $u\equiv0$.
Hence $(v,\te)$ satisfies
\begin{equation}\label{E43}
\begin{cases}
v_t-\Delta v=\nabla\te,\\
\te_t=\dv v,\\
v\vert_{t=0}=0,\quad \te\vert_{t=0}=\te_0.
\end{cases}
\end{equation}

Taking the Fourier transform $\sF$ in \eqref{E43}, we obtain
\[
\partial_t(\sF v)=-|\xi|^2\sF v+i\xi\sF\te,
\qquad
\partial_t(\sF\te)=i\xi\cdot\sF v,
\]
and thus,
\[
\begin{aligned}
\partial_t\bigl(|\sF v|^2+|\sF\te|^2\bigr)
=-2|\xi|^2|\sF v|^2
&\leq0,\\
\partial_t\left[e^{2t|\xi|^2}
\bigl(|\sF v|^2+|\sF\te|^2\bigr)\right]
=2|\xi|^2e^{2t|\xi|^2}|\sF\te|^2
&\ge0.
\end{aligned}
\]
Integrating from $0$ to $t$ and using $v(0)=0$ yields
\begin{equation}\label{E41}
e^{-2t|\xi|^2}|\sF\te_0(\xi)|^2
\le |\sF v(t,\xi)|^2+|\sF\te(t,\xi)|^2
\le |\sF\te_0(\xi)|^2.
\end{equation}
Fix $s \ge 1$. 
Integration of the $\dot H^s$ energy identity yields that
\begin{equation*}
    \frac12\|\Lambda^s(v,\te)(T)\|_{L^2}^2
+\int_t^T\|\Lambda^{s+1}v(\tau)\|_{L^2}^2\dtau=    \frac12\|\Lambda^s(v,\te)(t)\|_{L^2}^2.
\end{equation*}

Hence, for initial data satisfying the smallness condition in Theorem~\ref{Thm1}, the above identity shows that
\begin{equation}\label{E42}
 \frac12\|\Lambda^s(v,\te)(t)\|_{L^2}^2 =
\lim_{T\to\infty}\left(\frac12\|\Lambda^s(v,\te)(T)\|_{L^2}^2+ \int_t^T\|\Lambda^{s+1}v(\tau)\|_{L^2}^2\,\ud\tau \right)= \int_t^\infty\|\Lambda^{s+1}v(\tau)\|_{L^2}^2\,\ud\tau.
\end{equation}

Let $(v_n,\te_n)$ solve \eqref{E43} with data $(0,\psi_n)$. Since \eqref{E43} evolves each Fourier frequency independently,
\[
\sF v(t,\cdot) 
=
\sum_{n \in \bN}\sF v_n(t,\cdot) (\cdot), 
\quad\text{and}\quad
\sF \te(t,\cdot) 
=
\sum_{n \in \bN}\sF \te_n(t,\cdot) (\cdot).
\]
Thanks to \eqref{E41}, both $\sF v_n(t,\cdot)$ and
$\sF \te_n(t,\cdot)$ are supported in $A_n$ for all $t \ge 0$.
Since $1<t_n|\xi|^2 < 4$ on $A_n$, \eqref{E41} yields
\[
\begin{aligned}
t_n^s \|\Lambda^s(v,\te)(t_n)\|_{L^2}^2
&\ge t_n^s \int_{A_n} |\xi|^{2s}
\left(
|\sF v_n(t_n,\xi)|^2
+ |\sF \te_n(t_n,\xi)|^2
\right) \, d\xi \\
&\ge t_n^s \int_{A_n} |\xi|^{2s}
e^{-2t_n|\xi|^2} |\sF \psi_n(\xi)|^2 \, d\xi \\
&\ge e^{-8} t_n^s \int_{A_n} |\xi|^{2s}
|\sF \psi_n(\xi)|^2 \, d\xi \\
&= e^{-8} \ep^2 4^{-n}
\int_{1<|\eta|<2} |\eta|^{2s}
|\sF \psi(\eta)|^2 \, d\eta \\
&\ge e^{-8} \ep^2 4^{-n}
\int_{1<|\eta|<2} |\sF \psi(\eta)|^2 \, d\eta \\
&= e^{-8} \ep^2 4^{-n}.
\end{aligned}
\]
Hence, since $h(t_n) \le 4^{-n}$, we conclude that
\[
\frac{(1+t_n)^{s/2}}{h(t_n)}
\|\Lambda^s(v,\te)(t_n)\|_{L^2}
\ge e^{-4}\ep 2^n\longrightarrow\infty.
\]
Since neither $\psi$ nor the sequence $(t_n)$ depends on $s$,
the initial datum $\te_0$ constructed above satisfies this lower
bound for every real $s\ge1$.
Together with $u\equiv0$ and \eqref{E42}, this proves
Theorem~\ref{Thm2:sharp}.
Moreover, \eqref{E42} gives the same sharpness for the
dissipation tail alone.
\end{proof}

\subsection*{\textbf{Data Availability}}
We do not analyse or generate any datasets, because our work proceeds within a theoretical and mathematical approach.

\subsection*{\textbf{Conflict of Interest}}
All authors have no conflicts of interest.

\appendix 
\section{Time continuity of the temperature.}
\label{app:time_continuity}
Recall that the strong solution $(u, v, \te)$ to the system \eqref{E11} enjoys the following regularity:
\begin{align*}
(u, v) 
&\in \cC ([0, T];  H^1 (\bR^2)) \cap L^2 (0, T; H^2 (\bR^2)), \\
\theta 
&\in \cC ([0, T];  L^2 (\bR^2)) \cap L^\infty (0, T; H^1 (\bR^2)), \\
(\p_t u, \p_t v, \p_t \theta)
&\in L^2 (0, T; L^2 (\bR^2)), \\
\nabla u 
&\in L^1 (0, T; L^\infty (\bR^2)).
\end{align*}
In this section, we complete the proof of Theorem~\ref{T2}. That is, we show that $\theta \in \mathcal{C}([0,T];H^1(\bR^2))$. 
A standard argument using the stated regularity of $\te$ yields that $\te \in C_w([0,T];H^1(\bR^2))$.
The weak continuity and continuity of the norm imply the time continuity of the solution in a Hilbert space, and hence it remains only to prove the continuity of the map $t \mapsto \|\te(t)\|_{H^1}$.

We multiply $\eqref{E11}_2$-$\eqref{E11}_3$ by $(- \Delta) v$ and $(- \Delta) \theta$ respectively, and integrate over $\mathbb{R}^2$ to obtain
\begin{equation}
\label{appenn_A1}
\frac{1}{2} \frac{d}{dt}
\norm{\nabla (v, \theta)}_{L^2}^2
+ \norm{\Delta v}_{L^2}^2
= - \int_{\bR^2}
(u \cdot \nabla) v \cdot (- \Delta) v dx
- \int_{\bR^2}
(v \cdot \nabla) u \cdot (- \Delta) v dx
- \int_{\bR^2}
u \cdot \nabla \theta
(- \Delta) \theta dx.
\end{equation}
Using the divergence-free condition of $u$, integration by parts, and Young's inequality,
\begin{align*}
&\left| - \int_{\bR^2}
(u \cdot \nabla) v \cdot (- \Delta) v dx
- \int_{\bR^2}
(v \cdot \nabla) u \cdot (- \Delta) v dx
- \int_{\bR^2}
u \cdot \nabla \theta
(- \Delta) \theta dx \right| \\
&\le
C \norm{\nabla u}_{L^\infty}
\norm{\nabla v}_{L^2}^2
+ C
\norm{v}_{L^4}
\norm{\Delta u}_{L^2}
\norm{\nabla v}_{L^4}
+ C \norm{\nabla u}_{L^\infty}
\norm{\nabla \theta}_{L^2}^2\\
&\le
C \norm{\nabla u}_{L^\infty}
\norm{\nabla v}_{L^2}^2
+ C
\norm{v}_{L^2}^{\frac12}
\norm{\Delta u}_{L^2}
\norm{\nabla v}_{L^2}
\norm{\Delta v}_{L^2}^{\frac12}
+ C \norm{\nabla u}_{L^\infty}
\norm{\nabla \theta}_{L^2}^2\\
&\le
C \norm{\nabla u}_{L^\infty}
\norm{\nabla v}_{L^2}^2
+ C
\norm{v}_{L^2}^2
\norm{\nabla v}_{L^2}^4
+C\norm{\Delta u}_{L^2}^2
+\frac{1}{2}\norm{\Delta v}_{L^2}^2
+ C \norm{\nabla u}_{L^\infty}
\norm{\nabla \theta}_{L^2}^2.
\end{align*}

Applying \eqref{appenn_A1} to the fundamental theorem of calculus, we obtain
\begin{align*}
&\left| \norm{(v, \theta) (t)}_{\dot H^1}^2
- \norm{(v, \theta) (r)}_{\dot H^1}^2
+\int_r^t \norm{\Delta v (\tau)}_{L^2}^2 d \tau  \right| \\
&\le   
C\int_r^t \left(  \norm{\nabla u}_{L^\infty} \norm{\nabla v}_{L^2}^2
+ \norm{v}_{L^2}^2
\norm{\nabla v}_{L^2}^4
+\norm{\Delta u}_{L^2}^2
+ \norm{\nabla u}_{L^\infty} \norm{\nabla \theta}_{L^2}^2 \right) d \tau
,
\end{align*}
for any $0 \leq r \leq t < T$. The first and fourth terms on the right-hand side can be bounded
\begin{align*}
\int_r^t \norm{\nabla u}_{L^\infty} \left( \norm{\nabla v}_{L^2}^2 + \norm{\nabla \theta}_{L^2}^2 \right) d \tau
\le \norm{(v, \theta)}_{L^\infty(0, T; H^1)}^2 \int_r^t \norm{\nabla u}_{L^\infty} d \tau.
\end{align*}
For the second term and third term, it follows that
\begin{align*}
\int_r^t \left( \norm{v}_{L^2}^2
\norm{\nabla v}_{L^2}^4
+C\norm{\Delta u}_{L^2}^2\right) d \tau
&\le C \norm{v}_{L^\infty(0, T; H^1)}^6|t-r| + C\int_r^t \norm{\Delta u}_{L^2}^2  d \tau.
\end{align*}
Therefore, the entire right-hand side converges to $0$, which implies
\[
\left| \norm{(v, \theta) (t)}_{\dot H^1}^2 - \norm{(v, \theta) (r)}_{\dot H^1}^2 \right| \rightarrow 0 \quad\text{ as } r \to t.
\]
This proves that the map $t \mapsto \norm{(v, \theta) (t)}_{\dot H^1}$ is continuous on $[0, T]$, and therefore
\[
t \mapsto \norm{\theta (t)}_{H^1} \in C ([0, T]).
\]
The proof is finished.


\end{document}